 \documentclass[11pt]{article}
 \usepackage[top=1in,bottom=1in,left=1.5in,right=1.5in]{geometry}
  \usepackage{amsmath,amssymb}
  \usepackage{latexsym}
 \usepackage[dvips]{pstricks} 
  \usepackage{pst-node}
  \usepackage[dvips]{graphicx}

  \newcommand{\C}{\mathbb{C}}
  \newcommand{\F}{\mathbb{F}}
  \newcommand{\N}{\mathbb{N}}
  \renewcommand{\P}{\mathbb{P}}
  \newcommand{\R}{\mathbb{R}}
  \newcommand{\Z}{\mathbb{Z}}

  \newcommand{\e}{\mathbf{e}}
  \newcommand{\f}{\mathbf{f}}
  
  \newcommand{\gl}{\mathbf{GL}}

  \newcommand{\mP}{\mathbf{P}}
  \newcommand{\bP}{\mathbb{P}}

  \newcommand{\T}{\mathbf{T}}
  \newcommand{\U}{\mathbf{U}}
  \renewcommand{\u}{\mathbf{u}}
  \renewcommand{\v}{\mathbf{v}}

  \newcommand{\x}{\mathbf{x}}
  \newcommand{\X}{\mathbf{X}}
  \newcommand{\y}{\mathbf{y}}
  \newcommand{\z}{\mathbf{z}}
  \newcommand{\0}{\mathbf{0}}

  \newcommand{\cT}{\mathcal{T}}

  \newcommand{\rA}{\mathrm{A}}
  
  \newcommand{\rU}{\mathrm{U}}
  \newcommand{\rW}{\mathrm{W}}

  \newcommand{\hs}{\hspace*{\parindent}}
  \newcommand{\proof}{\hs \textbf{Proof.\ }}

  \newcommand{\Gr}{\mathop{\mathrm{Gr}}\nolimits}
  \newcommand{\trans}{^\top}
  \newcommand{\qed}{\hspace*{\fill} $\Box$\\}

  \newcommand{\inter}{\mathrm{int\;}}

  \newcommand{\brank}{\mathrm{brank\;}}
  \newcommand{\grank}{\mathrm{grank}}
  \newcommand{\mrank}{\mathrm{mrank}}
  \newcommand{\rmrank}{\mathrm{mtrank}}

  \newcommand{\proj}{\mathrm{proj\;}}
  
  \newcommand{\range}{\mathrm{range\;}}
  \newcommand{\Sing}{\mathrm{Sing\;}}

  \newcommand{\rD}{\mathrm{D}}

  \newcommand{\rS}{\mathrm{S}}

  \newcommand{\rank}{\mathrm{rank\;}}

  \newtheorem{theo}{\bfseries \hs Theorem}[section]
  \newtheorem{defn}[theo]{\bfseries \hs Definition}
  \newtheorem{prop}[theo]{\bfseries \hs Proposition}

  \newtheorem{corol}[theo]{\bfseries \hs Corollary}
  \newtheorem{con}[theo]{\bfseries \hs Conjecture}

  \numberwithin{equation}{section} 

 \renewcommand{\span}{\mathrm{span}}

\begin{document}

 \title{On the generic and typical ranks of 3-tensors}

 \author
 {Shmuel Friedland
 \thanks{Department of Mathematics, Statistics and Computer Science,
 University of Illinois at Chicago, Chicago, Illinois 60607-7045, USA,
 \emph{email}:friedlan@uic.edu} \thanks{Dedicated to the memory of Gene
 Golub}}
 \date{January 21, 2011}
 \maketitle

 \begin{abstract}
 We study the generic and typical ranks of $3$-tensors of dimension $l\times m \times n$ using results
 from matrices and algebraic geometry.  We state a conjecture
 about the exact values of the generic rank of $3$-tensors over
 the complex numbers, which is verified numerically for $l,m,n\le 14$.  We also discuss the typical ranks
 over the real numbers, and give an example of an infinite family of $3$-tensors of the form
 $l=m, n=(m-1)^2+1, m=3,4,\ldots$, which have at least two typical ranks.

 \end{abstract}


 \section{Introduction}

 The subject of tensors, their rank and the approximation of
 tensors by low rank tensors became recently a very active
 area of pure and applied mathematics.  See the reference section
 of this paper.
 $2$-dimensional tensors, which are identified as matrices,
 are well understood theoretically and numerically.  Tensors of dimension
 greater than $2$, are much more complicated theoretically and
 numerically than matrices.  Basically, matrices are strongly
 connected to linear operators, while tensors are strongly
 connected to the study of polynomial equations in several
 variables, which are best dealt with the tools of
 algebraic geometry.  Indeed, there is a vast literature in
 algebraic geometry discussing tensors.  See for example \cite[Chapter 20]{BCS}
 and references therein.
 Unfortunately, it is unaccessible to most researchers in applied and numerical
 analysis.

 The object of this paper three-fold. First, we give a basic
 introduction to one of the most interesting topics: the rank of
 $3$-tensors.  Second, we state our conjecture for the generic tensors of 3-tensors
 over the complex numbers. Third, we give general results for the typical ranks
 of 3-tensors over the real numbers.
 We illustrate the strength and generality of our results by comparing them to the known results
 in the literature.
 The novelle results of this paper are obtained by using results on matrices and basic results
 of algebraic geometry on polynomial equations over complex and real numbers.
 For reader's benefit we added a short appendix on complex and real algebraic geometry.
 The exact references for the results in complex and real algebraic geometry used in this paper are given in the appendix.

 This paper is written for the audience who has the knowledge
 of matrix theory and was only occasionally exposed to the
 study of polynomial maps in several complex variables.
 This paper is an expanded version of the talk I gave in
 \emph{Workshop on Algorithms for Modern Massive Data Sets}, sponsored by Computer
 Forum of the Stanford Computer Science Department, NSF and Yahoo!
 Research, June 21-24, 2006, \cite{Fr1}.

 We now survey briefly the contents of this paper.
 \S2 deals with the basic notions of the tensor product of
 three vector spaces over any field $\F$, $3$-tensors and their
 rank.  Theorem \ref{chartrank} gives a simple and useful characterization
 of the rank of a given tensor over any field $\F$, in terms of
 the minimal dimension of a subspace spanned by rank one
 matrices, containing a given subspace $L$ of $\F^{m\times n}$.
 \S3 introduces the notion of the generic rank in $\C^{l\times m
 \times n}$, denoted by $\grank(l,m,n)$.  ($\grank(l,m,n)$ is a
 symmetric function in $l,m,n$.)
  \S4 introduces the notion of the maximal rank in
 $\C^{l\times m \times n}$, denoted by $\mrank(l,m,n)$.
 \S5 gives known values for $\grank(l,m,n)$ and states the conjectured
 values of $\grank(l,m,n)$ in  the range $3\le l\le
 m\le n\le (l-1)(m-1)-1$.  This conjecture is verified numerically for all values of $l,m,n\le 14$.
 (Compare these results with the numerical results for $\grank(l,m,n)$ given in \cite[Table 1]{CBLC},
 for the values $l\le 4, m\le 5, n\le 12$.)
 \S6 shows how to apply some results
 on matrices to obtain bounds on $\grank(l,m,n)$ and
 $\mrank(l,m,n)$.  \S7 discusses the notion of typical ranks of
 real tensors $\R^{l\times m\times n}$, which are the analogs of generic rank over the complex numbers.  In this case one has a finite number of typical ranks taking all the values from $\grank(l,m,n)$
 to $\rmrank(l,m,n)$.  The typical ranks for the case $l=2\le m\le n$ are known.  For $m<n$ there is one typical rank which is equal to $\grank(2,m,n)=\min(n,2m)$.  For $2\le m=n$ there are
 two typical ranks $\grank(2,m,m)=m$ and $\rmrank(2,m,m)=m+1$.  See \cite{tB91} and \cite{tBK99}.
 In this paper we give another countable set of examples of the
 form $3\le l=m, n=(m-1)^2+1, m=3,\ldots$, where the maximal typical rank is strictly
 bigger than $\grank(m,m,(m-1)^2+1)=(m-1)^2+1$, i.e. there are at least two typical ranks in these cases.  The case $m=3$ is studied in \cite{tB04}.
 It is shown there that $\rmrank(3,3,5)=6$.  (It is not known that if $\rmrank(l,m,n)\le \grank(l,m,n)+1$,
 which holds in all known examples.)
 \S\ref{Sec:Appendix} gives a concise exposition of facts in complex and algebraic geometry needed here,
 with suitable references.

 \section{Basic notions and preliminary results}

 In this section we let $\F$ be any field.
  Usually we denote by a bold capital letter a finite dimensional vector
 space $\U$ over $\F$, unless stated otherwise.
 A vector $\u\in\U$ is denoted
 by a bold face lower case letter.  A matrix
 $A\in \F^{m\times n}$ denoted by a capital letter $A$, and we let either
 $A=[a_{ij}]_{i=j=1}^{m\times n}$ or simply $A=[a_{ij}]$.
 A \emph{3-tensor array}
 ${\cT}\in \F^{l\times m\times n}$ is
 denoted by a capital calligraphic letter.  So either
 ${\cT}=[t_{ijk}]_{i=j=k=1}^{l,m,n}$ or simply ${\cT}=[t_{ijk}]$.

 Let $\U_1,\U_2,\U_3$ be three vectors spaces.  Let $m_i:=\dim \U_i$ be the
 dimension of the vector space $\U_i$.  Let $\u_{1,i},
 \ldots,\u_{m_i,i}$ be a basis of $\U_i$ for $i=1,2,3$.
 Then $\U:=\U_1\otimes\U_2\otimes\U_3$ is the tensor product
 of $\U_1,\U_2,\U_3$.  $\U$ is a vector space of dimension
 $m_1m_2m_3$, and
 \begin{equation}\label{tenbas}
 \u_{i_1,1}\otimes\u_{i_2,2}\otimes\u_{i_3,3},\quad
 i_j= 1,\ldots,m_j, j=1,2,3,
 \end{equation}
 is a basis of $\U$.  For any
 permutation $\sigma:\{1,2,3\}\to\{1,2,3\}$ the tensor product
 $\U_{\sigma(1)}\otimes\U_ {\sigma(2)}\otimes \U_{\sigma(3)}$
 is isomorphic to $\U$.  Hence it will be convenient to assume
 that
 \begin{equation}\label{massump}
 1\le m_1\le m_2 \le m_3,
 \end{equation}
 unless stated otherwise.
 A \emph{3-tensor} is a vector in $\U$.
 We will call 3-tensor a  \emph{tensor}, and denote it by a Greek
 letter.  A tensor $\tau$ has the representation
 \begin{equation}\label{tenrep}
 \tau=\sum_{i_1=i_2=i_3=1}^{m_1,m_2,m_3}
 t_{i_1i_2i_3}\u_{i_1,1}\otimes\u_{i_2,2}\otimes \u_{i_3,3},
 \end{equation}
 in the basis (\ref{tenbas}).  If the basis (\ref{tenbas}) is fixed
 then $\tau$ is identified with ${\cT}=[t_{i_1i_2i_3}]\in
 \F^{m_1\times m_2\times m_3}$.

 Recall that $\x_1\otimes\x_2\otimes \x_3$, were $\x_i\in\U_i,
 i=1,2,3$, is called a \emph{rank one} tensor, or an
 \emph{indecomposable tensor}.  (Usually one assumes that all
 $\x_i\ne \0$.  Otherwise $\0=\x_1\otimes\x_2\otimes\x_3$ is called
 a rank zero tensor.)
 (\ref{tenrep}) is a decomposition
 of $\tau$ as a sum of at most $m_1m_2m_3$  rank one tensors, as
 $t_{i_1i_2i_3}\u_{i_1,1}\otimes\u_{i_2,2}\otimes \u_{i_3,3}=
 (t_{i_1i_2i_3}\u_{i_1,1})\otimes\u_{i_2,2}\otimes \u_{i_3,3}$.
 A decomposition of $\tau\ne\0$ to a sum of rank one tensors is given
 by
 \begin{equation}\label{rankonedec}
 \tau=\sum_{i=1}^k \x_i\otimes\y_i\otimes \z_i,\quad
 \x_i\in\U_1,\y_i\in\U_2,\z_i\in\U_3, \;i=1,\ldots,k.
 \end{equation}
 The minimal $k$ for which the above equality holds is called
 the \emph{rank} of the tensor $\tau$.  It is completely analogous
 to the rank of matrix $A=[a_{i_1i_2}]\in \F^{m_1\times m_2}$, which can be
 identified with 2-tensor in $\sum_{i_1=i_2=1}^{m_1,m_2} a_{i_1i_2} \u_{i_1,1}\otimes
 \u_{i_2,2}\in\U_1\otimes \U_2$.
 It is well known that, unlike in the case of matrices, the rank of a
 tensor may depend on the ground field $\F$.  In particular, by
 considering the algebraic closed field $\C$ versus $\R$, one may
 decrease the rank of the real valued tensor $\tau$.

 For $j\in\{1,2,3\}$  denote by
 $j^c:=\{p,q\}=\{1,2,3\}\backslash\{j\}$, where $1\le p <q\le 3$.
 Denote by $\U_{j^c}=\U_{\{p,q\}}:=\U_p\otimes\U_q$.
 A tensor $\tau\in \U_1\otimes\U_2\otimes \U_3$ induces a linear
 transformation $\tau(j):\U_{j^c}\to \U_j$ as follows.  Assume
 that $\u_{1,l},\ldots,\u_{m_l,l}$ is a basis in $\U_l$ for
 $l=1,2,3$.  Then any $\v\in \U_{j^c}$ is of the form
 $\v=\sum_{i_p=i_q=1}^{m_p,m_q} v_{i_pi_q} \u_{i_p,p}\otimes
 \u_{i_q,q}$.  Define
 \begin{equation}\label{deftauj}
 \tau(j)\;\v=\sum_{i_j=1}^{m_j}\big(\sum_{i_p,i_q=1}^{m_p,m_q} t_{i_1i_2i_3}
 v_{i_pi_q}\big)\u_{i_j,j}.
 \end{equation}
 The ${\mathrm{rank_j}} \tau$ \emph{is the rank of the operator} $\tau(j)$.
 Equivalently, let $A(j)=[a_{li_j}]\in \F^{m_pm_q\times m_j}$, where
 each integer $l\in [1,m_pm_q]$ corresponds to the pair $(i_p,i_q)$,
 for $ i_p=1,\ldots,m_p, i_q=1,\ldots,m_q$, and $i_j\in [1,m_j]\cap \N$.
 (For example arrange
 the pairs $(i_p,i_q)$ in the lexicographical order.
 Then $i_p=\lceil \frac{l}{m_q}\rceil$
 and $i_q=l-(i_p-1)m_q$.)  Set $a_{li_j}=t_{i_1i_2i_3}$.  Then
 ${\mathrm{rank_j}}\emph{}\tau=\rank A(j)$.  Associating a matrix $A(j)$ with
 the $3$-tensors is called \emph{unfolding} $\tau$ in direction
 $j$.  The following proposition is straightforward.
 \begin{prop}\label{rankdef}
 Let $\tau\in \U_1\otimes\U_2\otimes\U_3$ be given by
 (\ref{tenrep}).  Fix $j\in\{1,2,3\}, j^c=\{p,q\}$.  Let
 $T_{i_j,j}:=[t_{i_1i_2i_3}]_{i_p=i_q=1}^{m_p,m_q}\in \F^{m_p\times
 m_q}, i_j=1,\ldots,m_j$.  Then ${\mathrm{rank_j}}\tau$ is the dimension of subspace
 of $m_p\times m_q$ matrices spanned by $T_{1,j},\ldots,T_{m_j,j}$.
 \end{prop}
 The following result is well known.
 \begin{prop}\label{lowupbdrank}
 Let $\tau\in \U_1\otimes\U_2\otimes\U_3$.  Let
 $r_j:={\mathrm{rank_j}}\tau$ for $j=1,2,3$.  Denote by $0\le R_1\le R_2\le R_3$
 the rearranged values of $r_1,r_2,r_3$.  Then
 $R_3\le \rank \tau\le R_1R_2$.

 \end{prop}
 \proof  We first show that $r_3\le \rank \tau$.  Since
 $\u_{i_1,1}\otimes \u_{i_2,2}\in \U_{\{1,2\}}$ it follows
 that the decomposition (\ref{tenrep}) is a decomposition of
 $\tau_3$ to a sum of rank one linear operators from $\U_{\{1,2\}}$ to $\U_3$.
 Hence $r_3\le \rank \tau$.
 Let $j\in\{1,2,3\}, j^c=\{p,q\}$.  Recall that
 $\U$ is isomorphic to $\U':=\U_p\otimes \U_q\otimes \U_j$.
 Hence $r_j\le \rank \tau$ for $j=1,2$.  Thus $R_3\le \rank
 \tau$.

 Let $\v_{1,j},\ldots,\v_{1,r_j}$ be the basis of $\X_j:=\tau_j
 (\U_p\otimes\U_q)\subseteq \U_j$.  It is straightforward to show
 that $\tau\in \X_1\otimes\X_2\otimes \X_3$.  So $\tau_j:\X_p\otimes
 \X_q\to \X_j$.  Assume that $R_1=r_j$.  Decompose $\tau_j=\sum_{l=1}^{R_1} \z_l\otimes
 \x_l$, where $\z_l\in \X_p\otimes\X_q, \x_l\in \X_j$ for
 $l=1,\ldots, R_1$.  Since $\z_l\in \X_p\otimes\X_q$, it follows that
 each $\z_l$ is at most a sum of $R_2$ rank one tensors in
 $\X_p\otimes \X_q$.  Hence $\tau$ is a sum of at most $R_1R_2$
 rank one tensors in $\X_1\otimes\X_2\otimes\X_3$.  \qed \\
 The following proposition is obtained straightforward:
 \begin{prop}\label{changebas3t}
 Let the assumptions and the notations of Propositions
 \ref{rankdef}-\ref{lowupbdrank} hold.
 Let $[\v_{1,1},\ldots,\v_{m_1,1}],
 [\v_{1,2},\ldots,\v_{m_2,2}]$ be two bases in $\U_1,\U_2$
 respectively, where
 \begin{eqnarray*}
 &&[\u_{1,1},\ldots,\u_{m_1,1}]=[\v_{1,1},\ldots,\v_{m_1,1}]Q_1,\quad
 [\u_{1,2},\ldots,\u_{m_2,2}]=[\v_{1,2},\ldots,\v_{m_2,2}]Q_2,\\
 &&Q_1=[q_{pq,1}]_{p,q=1}^{m_1}\in\gl(m_1,\F),\quad
 Q_2=[q_{pq,2}]_{p,q=1}^{m_2}\in\gl(m_2,\F).
 \end{eqnarray*}
 Let
 $$\tau=\sum_{i,j,k=1}^{m_1,m_2,m_3} \tilde
 t_{ijk}\v_{i,1}\otimes \v_{j,2}\otimes\u_{j,3},\; \tilde
 T_{k,3}:=[\tilde t_{ijk}]_{i,j=1}^{m_1,m_2}\in \F^{m_1\times m_2},\;
 k=1,\ldots,m_3.$$
 Then $\tilde T_{k,3}=Q_1 T_{k,3} Q_2^T$ for
 $k=1,\ldots,m_3$.

 Let $[\v_{1,3},\ldots,\v_{m_3,3}]$ be another basis of $\U_3$,
 where
 $$[\u_{1,3},\ldots,\u_{m_3,3}]=[\v_{1,3},\ldots,\v_{m_3,3}]Q_3,\quad
 Q_3=[q_{pq,3}]_{p,q=1}^{m_3}\in \gl(m_3,\F).$$
 Then
 $\tau=\sum_{i,j,k=1}^{m_1,m_2,m_3}t'_{ijk}\u_{i,1}\otimes\u_{j,2}\otimes\v_{k,3}$
 and
 $T'_{k,3}=[t'_{ijk}]_{i,j=1}^{m_1,m_2}=\sum_{l=1}^k q_{kl,3} T_l$.

 Let $[\v_{1,i},\ldots,\v_{m_i,i}]$ be a basis in
 $\U_i$ such that $\tau_i \U_{i^c} =\span (\v_{1,i},\ldots,\v_{r_i,i})$
 for $i=1,2,3$.  Then $\tau=\sum_{i=j=k}^{m_1,m_2,m_3}\hat
 t_{ijk}\v_{i,1}\otimes\v_{j,2}\otimes \v_{k,3}$ and $\hat
 T_{k,3}:=[\hat t_{ijk}]_{i=j=1}^{m_1,m_2}\in\F^{m_1\times m_2}$ for
 $k=1,\ldots,m_3$.  Then $\hat T_{k,3}=0$ for $k>r_3$ and $\hat
 T_{1,3},\ldots,\hat T_{r_3,3}$ are linearly independent.
 Furthermore, each $\hat T_{k,3}=S_k\oplus 0:=\left[\begin{array}{cc}S_k &0\\
 0&0\end{array}\right]$, where $S_k\in \F^{r_1\times r_2}$ for
 $k=1,\ldots,r_3$.  Moreover, the span of $\range S_1,\ldots,\range
 S_{r_ 3}$ and the span of $\range S_1\trans,\ldots,\range S_{r_ 3}\trans$
 are $\F^{r_1}$ and $\F^{r_2}$ respectively.

 \end{prop}
 The following result is a very useful characterization of the
 rank of $3$-tensor.
%
 \begin{theo}\label{chartrank} Let $\tau\in \U_1\otimes\U_2\otimes\U_3$ be given by
 (\ref{tenrep}).  Fix $j\in\{1,2,3\}, j^c=\{p,q\}$.  Let
 $T_{i_j,j}:=[t_{i_1i_2i_3}]_{i_p=i_q=1}^{m_p,m_q}\in \F^{m_p\times
 m_q}, i_j=1,\ldots,m_j$.  Then $\rank\tau$ is the minimal dimension of a subspace
 of $m_p\times m_q$ matrices spanned by rank one matrices, which
 contains the subspace spanned by $T_{1,j},\ldots,T_{m_j,j}$.

 \end{theo}
 \proof  It is enough to prove the Proposition for the case $j=3$.
 Proposition \ref{lowupbdrank} and its proof yields that it is
 enough to consider the case where $r_3=m_3$, i.e.
 $T_{1,3},\ldots,T_{m_3,3}$ are linearly independent.
 Let $r$ be the dimension of the minimal subspace
 of $m_1\times m_2$ matrices spanned by rank one matrices, which
 contains the subspace spanned by $T_{1,3},\ldots,T_{m_3,3}$.

 Suppose that equality (\ref{rankonedec}) holds.  Since $r_3=m_3$
 it follows that $\z_{1},\ldots,\z_{k}$ span $\U_3$.  Without
 loss of generality we may assume that $\z_{1},\ldots,\z_{m_3}$
 form a basis in $\U_3$.  For each $l>m_3$ rewrite each $\z_{l}$
 as al linear combination of $\z_1,\ldots,\z_{m_3}$.  Thus
 \begin{equation}\label{rewrtau}
 \z_{l}=\sum_{p=1}^{m_3}b_{lj}\z_j,\;l=m_3+1,\ldots,k,\;
 \tau=\sum_{j=1}^{m_3} (\x_j\otimes\y_j +\sum_{l=m_3+1}^k
 b_{lj} \x_l\otimes\y_l)\otimes \z_j.
 \end{equation}
 Hence
 \begin{equation}\label{tp3for}
 T_{j,3}=\x_j\y_j\trans +\sum_{l=m_3+1}^k
 b_{lj} \x_l\y_l\trans, \quad j=1,\ldots,m_3.
 \end{equation}
 In particular, the subspace spanned by
 $T_{1,3},\ldots,T_{m_3,3}$ is contained in the subspace spanned
 by $k$ rank one matrices $\x_1\y_1\trans,\ldots,\x_k\y_k\trans$.
 Therefore $r \le k$, hence $r\le \rank\tau$.

 Assume now that there exist $\x_i\in\F^{m_1},\y_i\in\F^{m_2},
 i=1,\ldots,k$ such that
 $T_{p,3}=\sum_{i=1}^k a_{pi}\x_i\y_i\trans$ for $p=1,\ldots,m_3$.
 View $\x_i\y_i\trans$ as $\x_i\otimes\y_i$.
 Then
 \begin{equation}\label{antauid}
 \tau=\sum_{p=1}^{m_3}(\sum_{i=1}^k a_{pi}\x_i\otimes\y_i)\otimes\z_p=
 \sum_{i=1}^k \x_i\otimes\y_i\otimes (\sum_{p=1}^{m_3} a_{pi}\z_{p}).
 \end{equation}
 Hence $k\ge \rank \tau$.
 So $\rank \tau=r$.  \qed

 \section{Generic rank}

 From now and $\F$ is either the field of complex
 numbers $\C$ or the field of real numbers $\R$, unless stated otherwise.
 We refer the reader to \S\ref{Sec:Appendix} for the notations and results in algebraic geometry
 used in the sequel.
 Let $\x_i\in\C^{m_i}, i=1,2,3$.  Then a rank one tensor
 $\x_1\otimes\x_2\otimes\x_3$ is a polynomial map
 $\f:\C^{m_1+m_2+m_3}\to \C^{m_1\times m_2\times m_3}\equiv
 \C^{m_1m_2m_3}$, i.e. $\f(\x_1,\x_2,\x_3):=\x_1\otimes\x_2\otimes\x_3$.
 Thus we identify a vector
 $\z=(z_1,\ldots,z_{m_1+m_2+m_3})\trans\in \C^{m_1+m_2+m_3}$ with
 $(\x_1\trans,\x_2\trans,\x_3\trans)\trans$, which is also denoted by $(\x_1,\x_2,\x_3)$,
 and a vector $\y\in
 \C^{m_1m_2m_3}$ with
 ${\cT}=[t_{i_1i_2i_3}]_{i_1=i_2=i_3}^{m_1,m_2,m_3}\in \C^{m_1\times m_2\times m_3}$.
 (Here we arrange the three indices of
 $[t_{i_1i_2i_3}]$ in the lexicographical order.)
 Then $\rD\f$, the Jacobian matrix of partial derivatives is given as
 \begin{equation}\label{Jfeexpres}
 \rD\f(\x_1,\x_2,\x_3)
 =[A_1(\x_2,\x_3)|A_2(\x_1,\x_3)|A_3(\x_1,\x_2)]\in
 \C^{m_1m_2m_3\times (m_1+m_2+m_3)},
 \end{equation}
 is viewed as a block matrix,
 where $A_i\in \C^{m_1m_2m_3\times m_i}$ for $i=1,2,3$.
 More precisely, let
 $$\e_{i_j,j}=(\delta_{1i_j}, \ldots,\delta_{m_ji_j})\trans,
 i_j=1,\ldots,m_j$$
 be the standard bases in $\C^{m_j}$ for $j=1,2,3$.
 Then
 \begin{eqnarray}
 &&A_1(\x_2,\x_3)=[\e_{1,1}\otimes\x_2\otimes\x_3| \cdots|
 \e_{m_1,1}\otimes\x_2\otimes\x_3]\in \C^{m_1m_2m_3\times m_1}, \nonumber\\
 &&A_2(\x_1,\x_3)=[\x_1\otimes\e_{1,2}\otimes\x_3| \cdots |\x_1\otimes\e_{m_2,2}\otimes\x_3]
 \in \C^{m_1m_2m_3\times m_2},\\
 &&A_3(\x_1,\x_2)=[\x_1\otimes\x_2\otimes\e_{1,3}| \cdots |\x_1\otimes\x_2\otimes\e_{m_3,3}]
 \in \C^{m_1m_2m_3\times m_3}. \nonumber
 \end{eqnarray}

 So the $p-th$ column of $A_1(\x_2,\x_3)$ is the tensor
 $\e_{p,1}\otimes\x_2\otimes\x_3$.  Similar statements holds
 for $A_2(\x_1,\x_3)$ and $A_3(\x_1,\x_2)$.
 \begin{prop}\label{genrankjf}  Let $\x_i\in \C^{m_i}, i=1,2,3$, and
 denote by $\f:\C^{m_1}\times\C^{m_2}\times\C^{m_3}\to \C^{m_1\times
 m_2\times m_3}$ the map
 $\f(\x_1,\x_2,\x_3):=\x_1\otimes\x_2\otimes\x_3$.
 Identify $\C^{m_1}\times\C^{m_2}\times\C^{m_3}, \C^{m_1\times
 m_2\times m_3}$ with $\C^{m_1+m_2+m_3},\C^{m_1m_2m_3}$
 respectively.    Then
 \begin{equation}\label{ubrnkjf}
 \rank \rD\f(\x_1,\x_2,\x_3) \le m_1+m_2+m_3 -2.
 \end{equation}
 Equality holds for any
 $\x_i\ne \0$ for $i=1,2,3$.
 \end{prop}
 \proof  Let $A_1(\x_2,\x_3), A_2(\x_1,\x_3), A_3(\x_1,\x_2)$ be
 defined as in (\ref{Jfeexpres}).
 Note that
 $$\sum_{i_1=1}^{m_1} x_{i_1,1}\e_{i_1,1}\otimes\x_2\otimes\x_3=
 \sum_{i_2=1}^{m_2} x_{i_2,2}\x_1\otimes\e_{i_2,2}\otimes\x_3=
 \sum_{i_3=1}^{m_3}
 x_{i_3,3}\x_1\otimes\x_2\otimes\e_{i_3,3}=\x_1\otimes\x_2\otimes\x_3$$
 That is, the columns of $A_1(\x_2,\x_3)$, $A_2(\x_1,\x_3)$ and $A_3(\x_1,\x_2)$
 all span the vector $\x_1\otimes\x_2\otimes\x_3$.
 Hence the inequality (\ref{ubrnkjf}) holds.

 Choose $\x_1=\e_{1,1},
 \x_2=\e_{2,1},\x_3=\e_{1,3}$.  Then in $\rD\f(\e_{1,1},\e_{1,2},\e_{1,3})$
 the column $\e_{1,1}\otimes\e_{1,2}\otimes\e_{1,3}$ appears three
 times.  After deleting two columns
 $\e_{1,1}\otimes\e_{1,2}\otimes\e_{1,3}$,
 we obtain $m_1+m_2+m_2-2$ linearly independent columns, i.e.
 $ \rank \rD\f(\e_{1,1},\e_{1,2},\e_{1,3})= m_1+m_2+m_3 -2$.
 If $\x_1,\x_2,\x_3\ne \0$, then each $\x_i$ can be extended to a basis
 in $C^{m_i}$.  Hence equality holds in (\ref{ubrnkjf}).
 \qed

 Let $k$ be a positive integer and consider the map
 $\f_k:(\C^{m_1}\times\C^{m_2}
 \times\C^{m_3})^k\to \C^{m_1\times m_2\times m_3}$ given by
 \begin{eqnarray}
 &&\f_k(\x_{1,1},\x_{1,2},\x_{1,3},\ldots,\x_{k,1},\x_{k,2},\x_{k,3})=\sum_{l=1}^k
 \f(\x_{l,1},\x_{l,2},\x_{l,3})=\sum_{l=1}^k
 \x_{l,1}\otimes\x_{l,2}\otimes\x_{l,3}, \nonumber\\
 &&\x_{l,j}\in \C^{m_j}, \;j=1,2,3, \; l=1,\ldots,k \label{deffk}.
 \end{eqnarray}

 In this paper the closure of a set $S\subset \F^n$, denoted by Closure $S$, is the closure
 in the standard topology of $\F^n$.
 Since $\f_k$ is a polynomial map it follows, (see Appendix \S\ref{Sec:AppendixC}).
 \begin{defn}\label{defrkXk}  Let $Y_k\subseteq \C^{m_1\times m_2\times m_3}$
 be the closure of
 $\f_k((\C^{m_1}\times\C^{m_2} \times\C^{m_3})^k)$.  Denote by
 $r(k,m_1,m_2,m_3)$ the dimension of the variety
 $Y_k$.  Let $U_k\subsetneq Y_k$ be the constructible algebraic subset
 of $Y_k$, of dimension $r(k,m_1,m_2,m_3)-1$ at most, possibly an empty set, such that
 $\f_k((\C^{m_1}\times\C^{m_2} \times\C^{m_3})^k)=Y_k\backslash
 U_k$.

 $\cT\in\C^{m_1\times m_2\times m_3}$ has a border rank $k$ if $\cT\in Y_k\backslash Y_{k-1}$, where
 $Y_0=\{0\}$.  The border rank of $\cT$ is denoted by  $\brank \cT$.  $\cT$ is called rank ill conditioned if $\brank \cT<\rank \cT$.
 \end{defn}
 Clearly, $r(k,m_1,m_2,m_3)$ is a nondecreasing sequence in
 $k\in\N$.  (See for more details the proof of Theorem
 \ref{genrank3ten} and Theorem \ref{maxrank}.)
 The notion of \emph{border rank} was introduced in \cite{BCLR}.
 \begin{prop}\label{descric}  The set of all ill conditioned tensors $\cT\in\C^{m_1\times m_2\times m_3}$ of border rank $k$ equals to $U_k\setminus Y_{k-1}$. This set is a constructible algebraic set of dimension $ r(k,m_1,m_2,m_3)-1$ at most.
 \end{prop}
 \proof
 Recall that $Y_k\setminus Y_{k-1}$ is the set of tensors of border rank $k$.  Hence
 $\f_k((\C^{m_1}\times\C^{m_2} \times\C^{m_3})^k)\backslash Y_{k-1}$ is the set of all tensor
 whose rank and border rank are $k$.  By definition $Y_k$ is a disjoint union of $\f_k((\C^{m_1}\times\C^{m_2} \times\C^{m_3})^k)$ and $U_k$.  Hence the set of all ill conditioned tensors of border rank $k$
 is $U_k\setminus Y_{k-1}$.  Since $U_k$ is a constructible algebraic subset of $Y_k$, where $\dim U_k<\dim Y_k$, and $Y_{k-1}$ is an algebraic set, it follows from the results in Appendix \S\ref{Sec:AppendixC} that $U_k\setminus Y_{k-1}$ is a constructible algebraic set of dimension $\dim U_k$
 at most.
 \qed

 See \cite{DL08} for related results on rank ill conditioned tensors.
 The following theorem is a version of what is called in literature
 \emph{Terracini's} Lemma \cite{Ter}.
 \begin{theo}\label{genrank3ten}  Let $m_1,m_2,m_3\ge 2$ be three
 positive integers.
 Assume that $\e_{i_j,j}=(\delta_{1i_j},\ldots,\delta_{m_ji_j})\trans
 \in \C^{m_j}, i_j=1,\ldots,m_j$ is the standard basis in
 $\C^{m_j}$ for $j=1,2,3$.  Let $\grank(m_1,m_2,m_3)$ be the smallest
 positive integer $k$
 satisfying the following property.  There exist $3k$ vectors
 $\x_{l,1}\in\C^{m_1},\x_{l,2}\in \C^{m_2},\x_{l,3}\in \C^{m_3},
 l=1,\ldots,k$ such that the following $k(m_1+m_2+m_3)$ tensors span
 $\C^{m_1\times m_2\times m_3}$:
 \begin{eqnarray}\label{spanten}
 \e_{i_1,1}\otimes \x_{l,2}\otimes\x_{l,3},\;
 \x_{l,1}\otimes\e_{i_2,2} \otimes\x_{l,3},\;
 \x_{l,1}\otimes\x_{l,2} \otimes \e_{i_3,3},\\
 i_j=1,\ldots,m_j,\;
 j=1,2,3,\; l=1,\ldots,k.\nonumber
 \end{eqnarray}

 Then there exist three algebraic sets $U\subsetneq V\subseteq W\subsetneq
 \C^{m_1\times m_3\times m_3} \equiv \C^{m_1m_2m_3}$ such that the
 following holds.
 \begin{enumerate}
 \item
 Any ${\cT}=[t_{i_1i_2i_3}]\in \C^{m_1\times m_2\times
 m_3}\backslash U$ has rank $\grank(m_1,m_2,m_3)$ at most.
 \item
 Any ${\cT}=[t_{i_1i_2i_3}]\in \C^{m_1\times m_2\times
 m_3}\backslash V$ has exactly rank $\grank(m_1,m_2,m_3)$.

 \item
 Let ${\cT}=[t_{i_1i_2i_3}]\in \C^{m_1\times m_2\times
 m_3}\backslash W$.  Then $\rank {\cT}=\grank(m_1,m_2,m_3)$.
 Furthermore  the set of all $3\grank(m_1,m_2,m_3)$
 vectors
 $$\x_{l,1}\in\C^{m_1},\x_{l,2}\in\C^{m_2}, \x_{l,3}\in\C^{m_3},
 l=1,\ldots, \grank(m_1,m_2,m_3)$$
 satisfying the equality

 \begin{equation}\label{3tensordecvar}
 {\cT}=\sum_{l=1}^{\grank(m_1,m_2,m_3)}
 \x_{l,1}\otimes\x_{l,2}\otimes\x_{l,3}
 \end{equation}
 is a union of $\deg\f_k$ of pairwise disjoint varieties  $T_i({\cT})\subsetneq
 (\C^{m_1}\times\C^{m_2}\times\C^{m_3})^{\grank(m_1,m_2,m_3)}$
 of dimension $(m_1+m_2+m_3)\grank(m_1,m_2,m_3)-m_1m_2m_3$ for $i=1,\ldots,\deg\f_k$.
 View each rank one tensor $ \x_{l,1}\otimes\x_{l,2}\otimes\x_{l,3}$
 as a point in $(\C\backslash\{0\})\times\C\P^{m_1-1}\times \C\P^{m_2-1}\times
 \C\P^{m_3-1}$.  Then
 the set of all $\grank(m_1,m_2,m_3)$
 rank one tensors
 $$(\x_{1,1}\otimes\x_{1,2}\otimes\x_{1,3},\ldots,
 \x_{\grank(m_1,m_2,m_3),1}\otimes\x_{\grank(m_1,m_2,m_3),2}
 \otimes\x_{\grank(m_1,m_2,m_3),3})$$
 in $(\C\backslash\{0\})\times\P\C^{m_1-1}\times \P\C^{m_2-1}\times \P\C^{m_3-1}$ satisfying
 (\ref{3tensordecvar}) is a disjoint union of $\deg\f_k$ varieties each of dimension
 $(m_1+m_2+m_3-2)\grank(m_1,m_2,m_3)-m_1m_2m_3$.

 \end{enumerate}

 \end{theo}

 \proof
 (\ref{Jfeexpres}) yields that
 \begin{eqnarray}\label{jfkfor}
 &&\rD\f_k(\x_{1,1},\ldots,\x_{k,3})=
 [A_1(\x_{1,2},\x_{1,3})|A_2(\x_{1,1},\x_{1,3})|A_3(\x_{1,1},\x_{1,2})|\ldots\\
 &&|A_1(\x_{k,2},\x_{k,3})|A_2(\x_{k,1},\x_{k,3})|A_3(\x_{k,1},\x_{k,2})].\nonumber
 \end{eqnarray}
 Moreover the column space of $\rD\f_k$ is spanned by the vectors
 (\ref{spanten}).  As in the proof of the Proposition
 \ref{genrankjf}, generically the rank of
 $\rD\f_k(\x_{1,1},\ldots,\x_{k,3})$ is equal to $r(k,m_1,m_2,m_3)$.
 (See Appendix \S\ref{Sec:AppendixC}, top of page 21, for the definition of the term \emph{generically}.)
 Thus, there exists a strict algebraic set $X_k\subsetneqq (\C^{m_1}\times\C^{m_2}
 \times\C^{m_3})^k$ such $\rank \rD\f_k(\x_{1,1},\ldots,\x_{k,3})=r(k,m_1,m_2,m_3)$
 for any $(\x_{1,1},\ldots,\x_{k,3})\in (\C^{m_1}\times\C^{m_2}
 \times\C^{m_3})^k\backslash X_k$ and
 $\rank \rD\f_k(\x_{1,1},\ldots,\x_{k,3})<r(k,m_1,m_2,m_3)$
 for any $(\x_{1,1},\ldots,\x_{k,3})\in  X_k$.

 Let $k=1$.  Then Proposition \ref{genrankjf} yields
 that generically $\rank \rD\f_1(\x_{1,1},\x_{1,2},\x_{1,3})=m_1+m_2+m_3 -2$.
 Hence $f_1(\C^{m_1}\times\C^{m_2} \times\C^{m_3})$ is a
 constructible algebraic set of dimension $m_1+m_2+m_3-2$.  (In this case it is
 straightforward to show that $f_1(\C^{m_1}\times\C^{m_2} \times\C^{m_3})$
 is a variety.)  If $m_1+m_2+m_3-2=m_1m_2m_3$ then $f_1(\C^{m_1}\times\C^{m_2} \times\C^{m_3})=
 \C^{m_1\times m_2\times m_3}$, $\grank(m_1,m_2,m_3)=1$ and the
 theorem is trivial in this case.  That is every tensor $\cT$ is
 either rank one or rank zero tensor.

 Assume now that
 $m_1m_2m_3 > m_1+m_2+m_3-2$.  Then $f_1(\C^{m_1}\times\C^{m_2} \times\C^{m_3})\subsetneq
 \C^{m_1\times m_2\times m_3}$
 is a strict subvariety of tensors of rank 1 at most.
 Since
 $\f_k(\x_{1,1},
 \ldots,\x_{k,3})=\f_{k+1}(\x_{1,1},\ldots,\x_{k,3},\0,\0,\0)$,
 it follows
 \begin{eqnarray}\label{fkcon}
 \f_k((\C^{m_1}\times\C^{m_2} \times\C^{m_3})^k)\subseteq
 \f_{k+1}((\C^{m_1}\times\C^{m_2} \times\C^{m_3})^{k+1}),\;
 k=1,\ldots\\
 \textrm{ and } \f_k((\C^{m_1}\times\C^{m_2}
 \times\C^{m_3})^k)=\C^{m_1\times m_2\times m_3}
 \textrm{ for } k \ge m_1m_2m_3. \nonumber
 \end{eqnarray}
 In particular
 \begin{eqnarray}
 &&r(k,m_1,m_2,m_3), k=1,\ldots \textrm{ a
 nondecreasing sequence,}\nonumber\\
 &&r(\grank(m_1,m_2,m_3)-1, m_1,m_2,m_3)<m_1m_2m_3, \label{rkprop}\\
 &&r(k, m_1,m_2,m_3)=m_1m_2m_3 \textrm{ for } k \ge \grank(m_1,m_2,m_3).
 \nonumber
 \end{eqnarray}
 So
 $1<\grank(m_1,m_2,m_3)\le m_1m_2m_3$.
 Furthermore, $Y_{\grank(m_1,m_2,m_3)-1}$ is a strict subvariety
 of $\C^{m_1\times m_2\times m_3}$.
 Since $\C^{m_1\times m_2\times m_3}$ is the only variety of dimension
 $m_1m_2m_3$ in $\C^{m_1\times m_2\times m_3}$ it follows that
 $Y_k=\C^{m_1\times m_2\times m_3}$ for $k\ge \grank(m_1,m_2,m_3)$.

 Let $U:=U_{\grank(m_1,m_2,m_3)}$ as defined in Definition \ref{defrkXk}.
 Then any ${\cT}\in \C^{m_1\times m_2\times
 m_3}\backslash U$ is equal to some
 $\f_{\grank(m_1,m_2,m_3)}(\x_{1,1},\ldots,\x_{\grank(m_1,m_2,m_3),3})$,
 i.e. $\cT$ is of rank

 \noindent
 $\grank(m_1,m_2,m_3)$ at most.  This proves \emph{1}.

 Let $V=U\cup Y_{\grank(m_1,m_2,m_3)-1}$.  Then ${\cT}\in \C^{m_1\times m_2\times
 m_3}\backslash V$ has rank $\grank(m_1,m_2,m_3)$, i.e. \emph{2} holds.
 Let ${\cT}\in \C^{m_1\times m_2\times m_3}\backslash V$.  Then
 $\f_{\grank(m_1,m_2,m_3)}^{-1}(\cT)$ is a nonempty algebraic set
 of $(\C^{m_1}\times\C^{m_2} \times\C^{m_3})^{\grank(m_1,m_2,m_3)}$.
 As stated in \S\ref{Sec:AppendixC}, there exists a strict algebraic subset
 $W\subset \C^{m_1\times m_2 \times m_3}$, which contains $V$,
 such that the first claim of \emph{3} holds.

 Recall that rank one tensor
 $\x_{l,1}\otimes\x_{l,2}\otimes\x_{l,3}$ is a point in the manifold
 $(\C\backslash\{0\})\times\P\C^{m_1-1}\times \P\C^{m_2-1}\times \P\C^{m_3-1}$
 of dimension $m_1+m_2+m_3-2$.  Hence $\f_k$ can be viewed as a map
 $\tilde \f_k:((\C\backslash\{0\})\times\P\C^{m_1}\times \P\C^{m_2}\times
 \P\C^{m_3})^k\to \C^{m_1\times m_2 \times m_3}$.
 This interpretation of $\f_k$, combined with the first part of
 \emph{3} yields the second part of \emph{3}.  \qed
 \begin{defn}\label{defgenrank}
 $\;$
 \begin{itemize}
 \item The integer $\grank(m_1,m_2,m_3)$
 is called the \emph{generic rank} of $\cT\in \C^{m_1\times m_2\times
 m_3}$.
 \item $k \;(\le \grank(m_1,m_2,m_3))$ is called small if there is
 a rank $k$ tensor $\cT$ of the form (\ref{rankonedec})
 such that the Jacobian matrix at $\cT$ has rank $k(m_1+m_2+m_3-2)$.
 \item $k \;(\ge \grank(m_1,m_2,m_3))$ is called big if there is a
 rank $k$ tensor $\cT$ of the form (\ref{rankonedec})
 such that the Jacobian matrix at $\cT$ has rank equal to the maximal rank $m_1m_2m_3$.
 \item $(m_1,m_2,m_3)$ is called perfect if
 $k=\frac{m_1m_2m_2}{m_1+m_2+m_3-2}$ is a small integer.
 \end{itemize}
 \end{defn}
 \begin{corol}\label{maxbordrank}   $\brank \cT\le \grank(m_1,m_2,m_3)$ for any $\cT\in\C^{m_1\times
 m_2\times m_3}$.
 \end{corol}

 The generic rank $\grank(m_1,m_2,m_3)$ has the following interpretation.
 Assume that the entries of $\cT\in\C^{m_1\times m_2\times m_3}$ are independent random variables,
 with normal complex Gaussian distribution.  Then with probability $1$ the rank of $\cT$
 is $\grank(m_1,m_2,m_3)$.  Furthermore, Proposition \ref{descric} yields that with probability
 $1$ the border rank of $\cT$ is also equal to $\grank(m_1,m_2,m_3)$.

 Since the dimension of any algebraic variety is nonnegative the second
 part of \emph{3} of Theorem \ref{genrank3ten}  yields the well known result, e.g.
 \cite[Chapter 20]{BCS}:
 \begin{corol}\label{grankub} $\grank(m_1,m_2,m_3)\ge
 \lceil\frac{m_1m_2m_3}{m_1+m_2+m_3-2}\rceil$.
 \end{corol}
 The following result is known, e.g. \cite[Prop. 2.3]{Str}, and we
 give its proof for completeness.
 \begin{prop}\label{incrgrank}  Let $m_1\ge l_1, m_2\ge l_2, m_3\ge l_3$ be positive
 integers. Then
 $\grank(m_1,m_2,m_3)\ge \grank(l_1,l_2,l_3)$.

 \end{prop}
 \proof  Since $\grank(m_1,m_2,m_3)$ is a symmetric function in
 $m_1,m_2,m_3$, it is enough to to show that $\grank(m_1,m_2,m_3),
 m_1=1,2,\ldots$ is a nondecreasing sequence.
 Assume that  $(T_{1,1},\ldots,T_{l+1,1})\in (\C^{m_2\times
 m_3})^{l+1}$ is a generic point.  Then $(T_{1,1},\ldots,T_{l,1})\in (\C^{m_2\times
 m_3})^{l}$ is also a generic point.  Theorem \ref{chartrank}
 implies that the minimal dimensions of subspaces spanned by rank
 one matrices containing $\span (T_{1,1},\ldots,T_{l+1,1}),\span (T_{1,1},\ldots,T_{l,1})$
 are $\grank(l+1,m_2,m_3)$, $\grank(l,m_2,m_3)$.   Hence $\grank(l+1,m_2,m_3)\ge
 \grank(l,m_2,m_3)$.   \qed

 \begin{prop}\label{tBSin}  Let $l\ge 3, m\ge 4$ be integers.  Then $\grank(l,m,m)\ge m+2$.
 \end{prop}
 \proof  Fix $m\ge 4$ and let $\phi(t)= \frac{tm^2}{t+2m-2}$ be a function of $t>0$.
 Then $\phi(t)$ is increasing.  Hence for $t\ge 3$
 \[\phi(t)\ge \phi(3)=\frac{3m^2}{2m+1}> m+1 \textrm{ for } m\ge 4.\]
 Therefore for $l\ge 3,m\ge 4$ $\grank(l,m,m)\ge m+2$. \qed

 As $\grank(3,3,3)=5$, see (\ref{(3,2p+1,2p+1)}) it follows that $\grank(l,m,m)\ge m+2$
 for $l,m\ge 3$, which was shown in \cite{tBS06}.

 \section{Maximal rank}\label{sec:mrank}
 \begin{theo}\label{maxrank}  Let $m_1,m_2,m_3,k$ be three positive
 integers and assume that $\f_k$ is given by (\ref{deffk}).
 Let $\mrank(m_1,m_2,m_3)$ be the smallest integer $k$ such that
 equality holds in (\ref{fkcon}).  I.e.
 \begin{equation}\label{maxrank1}
 \f_{k}((\C^{m_1}\times\C^{m_2} \times\C^{m_3})^k)=
 \f_{k+1}((\C^{m_1}\times\C^{m_2} \times\C^{m_3})^{k+1})
 \end{equation}
 for $k=\mrank(m_1,m_2,m_3)$, and
 \begin{equation}\label{maxrank2}
 \f_k((\C^{m_1}\times\C^{m_2} \times\C^{m_3})^k)\subsetneq
 \f_{k+1}((\C^{m_1}\times\C^{m_2} \times\C^{m_3})^{k+1}),
 \end{equation}
  for $k<\mrank(m_1,m_2,m_3)$.
 Then the maximal rank
 of all $3$-tensors in $\C^{m_1\times m_2\times m_3}$ is
 $\mrank(m_1,m_2,m_3)$, and
 \begin{equation}\label{maxrank3}
 \grank(m_1,m_2,m_3)\le
 \mrank(m_1,m_2,m_3).
 \end{equation}
 For each integer $k\in [1,\mrank(m_1,m_2,m_3)]$ the set of all
 tensors of rank $k$ is a nonempty constructible algebraic set
 $\f_k((\C^{m_1}\times\C^{m_2} \times\C^{m_3})^k)\backslash
 \f_{k-1}((\C^{m_1}\times\C^{m_2} \times\C^{m_3})^k)$, (
 $\f_0((\C^{m_1}\times\C^{m_2} \times\C^{m_3})^0):=\{\0\}$).
 If strict inequality in (\ref{maxrank3}) holds
 then the set of all $3$-tensors in
 $\C^{m_1\times m_2\times m_3}$ of rank greater than $\grank(m_1,m_2,m_3)$
 is a constructible algebraic set of $\C^{m_1\times m_2\times m_3}$ of dimension $m_1m_2m_3-1$ at most.
 Furthermore for each nonnegative integer
 $k<\grank(m_1,m_2,m_3)$ the following holds:
  \begin{equation}\label{dimin}
 \dim \f_k((\C^{m_1}\times\C^{m_2} \times\C^{m_3})^k) <
 \dim \f_{k+1}((\C^{m_1}\times\C^{m_2} \times\C^{m_3})^{k+1}).
 \end{equation}
 In particular for $k\le \grank(m_1,m_2,m_3)$ the dimension of the
 constructible algebraic set of all $3$-tensor of rank $k$ is
 \begin{equation}\label{dim=rk}
 \dim \f_k((\C^{m_1}\times\C^{m_2}
 \times\C^{m_3})^k)=r(k,m_1,m_2,m_3),
 \end{equation}
 which is the rank of the Jacobian matrix
 $\rD\f_k$ at the generic point $(\x_{1,1},\ldots,\x_{k,3})\in \C^{m_1\times m_2\times m_3}$,
 (which is also the maximal rank of
 $\rD\f_k(\x_{1,1},\ldots,\x_{k,3})$).

 \end{theo}

 \proof Assume the notation of Definition \ref{defrkXk} for
 $k\ge 0$, where $Y_0:=\{\0\}\in \C^{m_1\times m_2\times m_3}, U_0=\emptyset$.
 Suppose that (\ref{maxrank1}) holds for $k=p$.  Then any tensor of
 the form $\sum_{l=1}^{p+1} \x_{l,1}\otimes\x_{l,2}\otimes \x_{l,3}$
 is of the form $\sum_{l=1}^{p} \y_{l,1}\otimes\y_{l,2}\otimes
 \y_{l,3}$.  Hence the rank of any tensor is $p$ at most.  Thus
 (\ref{maxrank1}) holds for any $k\ge p$.  The second part of
 (\ref{fkcon}) yields $\mrank(m_1,m_2,m_3)\le m_1m_2m_3$,
 and $\f_k((\C^{m_1}\times\C^{m_2} \times\C^{m_3})^k)
 =\C^{m_1\times m_2\times m_3}$ for
 $k=\mrank(m_1,m_2,m_3)$.  Thus the rank of any $3$-tensor is at
 most $\mrank(m_1,m_2,m_3)$.  From the definition of
 $\mrank(m_1,m_2,m_3)$ we deduce (\ref{maxrank2}).
 That is for each integer $k\in [1,\mrank(m_1,m_2,m_3)]$,
 $Z_k:=(Y_k\backslash U_k)\backslash( Y_{k-1}\backslash U_{k-1})$ is
 the nonempty constructible algebraic set of rank $k$ tensors.

 From the definition of $q:=\grank(m_1,m_2,m_3)$ we deduce that
 $Y_k =\C^{m_1\times m_2\times m_3}$
 for $k\ge q$.  Hence $\f_k((\C^{m_1}\times\C^{m_2}
 \times\C^{m_3})^k)=\C^{m_1\times m_2 \times m_3}\backslash
 U_k, k\ge q$, where each $U_k$ for $k\ge q$ is a constructible algebraic
 set satisfying
 $$U_{q}\supsetneq U_{q+1}\supsetneq\ldots \supsetneq
 U_{\mrank(m_1,m_2,m_3)}=\emptyset.$$
 (Note that $U_k=\emptyset$ for $k>\mrank(m_1,m_2,m_3)$.)

 We now show (\ref{dimin}) for $k<q$.
 Definition \ref{defrkXk} implies the equality
 (\ref{dim=rk}).  Assume to the contrary that
 $r(k,m_1,m_2,m_3)=r(k+1,m_1,m_2,m_3)$ for some integer $k\in
 [1,q-1]$.  Let $s$ be the smallest positive integer satisfying
 this condition.  Then there exists an algebraic set
 $X_s\subsetneq (\C^{m_1}\times\C^{m_2} \times\C^{m_3})^s$ such that
 $\rank \rD\f_s(\x_{1,1},\ldots,\x_{s,3})=r(s,m_1,m_2,m_3)$
 for any $(\x_{1,1},\ldots,\x_{s,3})\in (\C^{m_1}\times\C^{m_2}
 \times\C^{m_3})^s\backslash X_s$.
 I.e., $s(m_1+m_2+m_3)$ tensors given in (\ref{spanten})
 span $r(s,m_1,m_2,m_3)$ dimensional subspace in $\C^{m_1\times m_2
 \times m_3}$ for any $(\x_{1,1},\ldots,\x_{s,3})\in (\C^{m_1}\times\C^{m_2}
 \times\C^{m_3})^s\backslash X_s$.

 Let
 $(\x_{1,1},\ldots,\x_{s+1,3})\in (\C^{m_1}\times\C^{m_2}
 \times\C^{m_3})^{s+1}$.  Then
 $$\rank \rD\f_{s+1}(\x_{1,1},\ldots,\x_{s+1,3})
 \le r(s+1,m_1,m_2,m_3)=r(s,m_1,m_2,m_3).$$
 I.e., $(s+1)(m_1+m_2+m_3)$ tensor given in (\ref{spanten})
 span at most $r(s,m_1,m_2,m_3)$ dimensional subspace in $\C^{m_1\times m_2
 \times m_3}$.
 Assume that $(\x_{1,1},\ldots,\x_{s,3})\in (\C^{m_1}\times\C^{m_2}
 \times\C^{m_3})^s\backslash X_s$.  Then  $(s+1)(m_1+m_2+m_3)$ tensor given in (\ref{spanten})
 span exactly $r(s,m_1,m_2,m_3)$ dimensional subspace in $\C^{m_1\times m_2
 \times m_3}$.  Moreover, the $s(m_1+m_2+m_3)$ tensor given by (\ref{spanten})
 for $k=s$ span the above subspace.  Hence the tensors
 $$\e_{i_1,1}\otimes \x_{s+1,2}\otimes\x_{s+1,3},
 \x_{s+1,1}\otimes\e_{i_2,2} \otimes\x_{s+1,3},
 \x_{s+1,1}\otimes\x_{s+1,2} \otimes \e_{i_3,3}, i_j=1,\ldots,m_j,
 j=1,2,3,$$
 are spanned by $s(m_1+m_2+m_3)$ tensor given by (\ref{spanten})
 for $k=s$.

 Let $k> s+1$ and consider $\rank \rD
 \f_k(\x_{1,1},\ldots,\x_{k,3})$, which is equal to the dimension
 of the subspace spanned by  $k(m_1+m_2+m_3)$ tensor given by
 (\ref{spanten}).  Let $(\x_{1,1},\ldots,\x_{s,3})\in (\C^{m_1}\times\C^{m_2}
 \times\C^{m_3})^s\backslash X_s$.  Then the above arguments show that
 $\rank \rD\f_{k}(\x_{1,1},\ldots,\x_{k,3})
 =r(s,m_1,m_2,m_3).$
 Since $X_s\times (\C^{m_1}\times\C^{m_2}
 \times\C^{m_3})^{k-s}$ is an algebraic set of $(\C^{m_1}\times\C^{m_2}
 \times\C^{m_3})^{k}$ it follows that $r(k,m_1,m_2,m_3)=r(s,m_1,m_2,m_3)$.
 This is impossible, since $r(s,m_1,m_2,m_3)<
 m_1m_2m_3=r(q,m_1,m_2,m_3)$.
 Hence (\ref{dimin}) holds for $k<q$.  \qed

 Combine the arguments of the proof of Theorem \ref{genrank3ten} with the results in \S\ref{Sec:AppendixC} to obtain.
 \begin{theo}\label{rankkrep}  Let $m_1,m_2,m_3,k$ be three positive
 integers and assume that $\f_k$ is given by (\ref{deffk}).
 Suppose that $k\le \grank(m_1,m_2,m_3)$.  Let $\cT\in \C^{m_1\times
 m_2\times m_3}$ be a generic tensor of rank $k$, i.e. a generic
 point in $\f_k((\C^{m_1}\times \C^{m_2}\times \C^{m_3})^k)\subset \C^{m_1\times m_2\times
 m_3}$.  Then the set of all possible decompositions of $\cT$ as a
 sum of $k$ rank one tensors is a disjoint union of $\deg f_k$ varieties of dimension
 $k(m_1+m_2+m_3 -2)-r(k,m_1,m_2,m_3)$.  In particular, if
 $r(k,m_1,m_2,m_3)=k(m_1+m_2+m_3 -2)$, i.e. $k$ is small, then $\cT$ can be
 decomposed as a sum of $k$-rank tensors in a finite number of ways
 given by a number $N(k,m_1,m_2,m_3)=\deg \f_k$.
 \end{theo}

 We remark that in the case $r(k,m_1,m_2,m_3)=k(m_1+m_2+m_3-2)$
 the positive integer
 $N(k,m_1,m_2,m_3)$ is divisible by $k!$, since we can permute
 the $k$ summands in (\ref{rankonedec}).  If $N(k,m_1,m_2,m_3)=k!$,
 this means that a generic rank $k$ tensor $\cT$ has a unique
 decomposition to $k$ factors.
 As we can see later,the numerical
 evidence points out that the equality $r(k,m_1,m_2,m_3)=k(m_1+m_2+m_3-2)$
 occurs for many
 $k<\grank(m_1,m_2,m_3)$.

 \section{Known theoretical results}\label{sec:knres}

 The following results are known.  See the references below.
 \begin{eqnarray}
 &&\grank(m_1,m_2,m_3)=\min (m_3,m_1m_2)  \textrm{ if } m_3\ge
 (m_1-1)(m_2-1)+1,\label{baseq}\\
 &&\textrm{in particular } \grank(2,m_2,m_3)=\min(m_3,2m_2) \textrm{ if } 2\le m_2\le m_3,
 \label{2case}\\
 &&\grank(3,2p,2p)=\lceil\frac{12p^2}{4p+1}\rceil \textrm{ and }
 \lfloor\frac{12p^2}{4p+1}\rfloor \textrm{ is small},
 \label{(3,2p,2p)}\\
 &&\grank(3,2p+1,2p+1)=\lceil\frac{3(2p+1)^2}{4p+3}\rceil+1,
 \label{(3,2p+1,2p+1)}\\
 &&(n,n,n+2) \textrm{ is perfect for } n\ne 2 \textrm{ (mod 3)},
 \label{(n,n,n+2)}\\
 &&(n-1,n,n) \textrm{ is perfect for } n=0 \textrm{ (mod 3)},
 \label{(n-1,n,n)}\\
 &&\grank(4,m,m)=\lceil\frac{4m^2}{2m+2}\rceil, \label{(4,n,n)}\\
 &&\grank(n,n,n)=\lceil\frac{n^3}{3n-2}\rceil
 \textrm{ and }
  \lfloor\frac{n^3}{3n-2}\rfloor
 \textrm{ is small for }  n\ne 3
 \label{(n,n,n)},\\
 &&(m_1,2m_2',2m_3') \textrm{ perfect if } \frac{2m_1m_2'}{m_1+2m_2'+2m_3'-2}
 \textrm{ is integer }, \label{perfect}\\
  &&\textrm{where } (\ref{massump}) \textrm{ holds}.\nonumber
 \end{eqnarray}

 See \cite{CGG} for (\ref{baseq}),  \cite{Str} for
 (\ref{(3,2p,2p)}- \ref{(n-1,n,n)}), \cite{BCS} for (\ref{(4,n,n)}),
 \cite{Lic85} and \cite[Theorem 5.3]{AOP} for (\ref{(n,n,n)}-\ref{perfect}).
 Note that in view of (\ref{baseq})
 \begin{equation}\label{stanper}
 (m_1,m_2, (m_1-1)(m_2-1)+1) \textrm{ is perfect}.
 \end{equation}
 We bring another proof of (\ref{baseq}) using matrices in
 \S\ref{sec:examples}.
 It was conjectured in \cite{Fr1}.
 \begin{con}\label{mcon} Let $3\le m_1\le m_2\le m_3\le (m_1-1)(m_2-1)$ and $(m_1,m_2,m_3)\ne
 (3,2p+1,2p+1), p\in\N$.  Then $\grank(m_1,m_2,m_3)=
 \lceil\frac{m_1m_2m_3}{m_1+m_2+m_3-2}\rceil$.

 \end{con}

 Combine Corollary \ref{grankub},  Proposition \ref{incrgrank}
 and (\ref{stanper}) to deduce.
 \begin{eqnarray}\label{(m,n,(m-1)(n-1))}
 \grank(m_1,m_2,m_3)=\lceil\frac{m_1m_2m_3}{m_1+m_2+m_3-2}\rceil=
 \grank(m_1,m_2,m_3+1)\\
 \textrm{ for } m_3=(m_1-1)(m_2-1) \textrm{ and }  3\le m_1,m_2.
 \nonumber
 \end{eqnarray}
 I.e., the above conjecture holds for $m_3=(m_1-1)(m_2-1)$.
 A more precise version of Conjecture \ref{mcon} is
 \begin{con}\label{kcon}.  Let the assumptions of Conjecture
 \ref{mcon} hold.  Then any integer $k\in
 [2,\lceil\frac{m_1m_2m_3}{m_1+m_2+m_3-2}\rceil-1]$ is small.
 \end{con}
  We call $(m_1,m_2,m_3)$ \emph{regular} if $(m_1,m_2,m_3)$ satisfies
 Conjecture \ref{mcon} and $\lfloor\frac{m_1m_2m_3}{m_1+m_2+m_3-2}\rfloor$
 is small.

 We verified numerically\footnote{I thank  M. Tamura for programming the software to compute
 $\grank(m_1,m_2,m_3)$  and $r(k,m_1,m_2,m_3)$.}
  the above two conjectures for $m_1\le m_2\le m_3\le
 14$ as follows. We chose at random $k\in [2,
 \lceil\frac{m_1m_2m_3}{m_1+m_2+m_3-2}\rceil]$ vectors $\x_{l,i}\in
 (\Z\cap [-99,99])^{m_i}, i=1,2,3, \;l=1,\ldots,k$ such that the
 rank of the Jacobian matrix at the corresponding rank $k$ tensor
 \begin{equation}\label{3tensordecvar1}
 {\cT}=\sum_{l=1}^{k}
 \x_{l,1}\otimes\x_{l,2}\otimes\x_{l,3}
 \end{equation}
 was
 $\min(k(m_1+m_2+m_3-2),m_1m_2m_3)$.  See also \cite{CBLC} for
 numerical results.

 The values of $\mrank(m_1,m_2,m_3)$ are much harder to compute.
 The following results are known.  First,  \cite[p'10]{Kr89}, (see also \cite{Ja79}),
 \begin{equation}\label{maxrnkpq2}
 \mrank (2,m,n)=m+\min(m, \lfloor \frac{n}{2}\rfloor) \textrm{ for } 2 \le m\le n.
 \end{equation}
 Second, it is claimed in \cite{Roc93} that
 \begin{equation}\label{maxrank333}
 \mrank(3,3,3)=5
 \end{equation}

 \section{Matrices and the rank of $3$-tensors}\label{sec:examples}

 In this section we use known results for matrices to
 find estimates on the generic and maximal rank of tensors.
 \begin{prop}\label{rankeqbm}  Let $\U_i$ be $m_i$-dimensional
 vector space over $\F$, for $i=1,2,3$.  Then
 $\mrank(m_1,m_2,m_3)=m_1m_2$ for $m_1m_2\le m_3$.
 More precisely,
 let $\tau\in \U_1\otimes\U_2\otimes\U_3$ be given by
 (\ref{tenrep}).  Let $R_1,R_2,R_3$ be defined as in Proposition
 \ref{lowupbdrank}.
 Assume that $R_3= R_1R_2$.  Then $\rank \tau = R_1R_2$.

 \end{prop}

 \proof  Since $\F^{m_1\times m_2}$ is spanned by $m_1m_2$ rank
 one matrices, Theorem \ref{chartrank} yields
 that $\mrank(m_1,m_2,m_3)\le m_1m_2$.
 Choose $\tau$ represented by (\ref{tenrep}), such that
 $T_{1,3},\ldots,T_{m_3,3}\in \F^{m_1\times m_2}$ span $\F^{m_1\times m_2}$.
 Theorem \ref{chartrank} yields that $\rank \tau =m_1m_2$, i.e.
 $\mrank(m_1,m_2,m_3) = m_1 m_2$.
 The second part of the proposition follows from Proposition
 \ref{lowupbdrank}.  \qed

 (The above results in this section hold for any field $\F$.  We remind the reader that from now and on $\F=\R,\C$.)
 We now show how to deduce (\ref{baseq}) using matrices.
 For a finite dimensional vector space $\U$ over $\F$ of dimension $N$ denote by $\Gr(k,\U)$,
 the $k$-Grassmannian,
 the manifold of all $k$ dimensional subspaces of $\U$.  ($k\in [0,N]$.)
 Note that $\Gr(1,\F^{m\times n})$ can be identified with $\P\F^{mn-1}$, a the projective space of dimension $mn-1$.  Equivalently, if $0_{m\times n}\ne A\in \F^{m\times n}$,
 then $\hat A\in \Gr(1,\F^{m\times n})$ corresponds to all points $tA,
 t\in \F\backslash \{0\}$.  Note that $\rank A =\rank tA$ for any
 $t\in\F\backslash\{0\}$.  Thus we define $\rank \hat A :=\rank A$.
 Usually we will identify $\hat A\in \Gr(1,\F^{m\times n})$ with one of $tA\in
 \F^{m\times n}\backslash \{0\}$ and no ambiguity will arise.

 Let $L \subseteq \F^{m\times n}$ be a
 subspace of dimension $d\ge 1$.  Then $\proj L\subset   \Gr(1,\F^{m\times
 n})$, the set of all one dimensional subspaces in $L$.  The
 dimension of $\proj L$ is $d-1$ and $\proj L$ can be identifies with
 $\P\F^{d-1}$.  $\proj L$ is called a \emph{linear space} in $\proj\F^{m\times
 n}$.
 The following result is known \cite{HT1, FK}.
 \begin{theo}\label{rnkvarmat}
 Let  $\rU_{k,m,n}(\F)\subseteq \F^{m\times n}$ be the set of all
 $m\times n$ matrices of rank $k$ at most.  Then $\rU_{k,m,n}(\F)$ is an
 irreducible variety of dimension $k(m+n-k)$.  Furthermore,
 $\rU_{k,m,n}(\F)\backslash \rU_{k-1,m,n}(\F)$ is quasi-projective variety
 of all matrices of rank $k$ exactly, which is a manifold of dimension
 $k(m+n-k)$.

 Any complex subspace of $L\subset \C^{m\times n}$ of dimension
 $(m-k)(n-k)+1$ contains a nonzero matrix of rank $k$ at most.
 More precisely, for a generic subspace
 $L\subset \C^{m\times n}$
 of dimension  $(m-k)(n-k)+1$, the linear space $\proj L$ contains exactly
 \begin{equation}
 \gamma_{k,m,n}:=
  \prod _{j=0} ^{n-k-1}\frac {\binom {m+j}{m-k}} {\binom {m-k+j}{m-k}}=
 \prod_{j=0}^{n-k-1}\frac{(m+j)!\,j!}{(k+j)!\,(m-k+j)!},
 \label{degkmn}
 \end{equation}
 distinct matrices of rank $k$ exactly.

 \end{theo}

 \begin{theo}\label{subspr1m}  Let $2\le m,n$ and $d\in
 [(m-1)(n-1)+1, mn-1]$ be fixed integers.  Then a generic
 subspace $L\subset \C^{m\times n}$ of dimension $d$
 is spanned by rank one matrices.
 \end{theo}
 \proof We first consider the case $d=(m-1)(n-1)+1$.  It is not
 difficult to check that $d\le \gamma_{1,m,n}$.
 Let $L$ be a generic subspace $L$ of dimension $(m-1)(n-1)+1$
 Then $L\cap U_{k,m,n}(\F)=\{A_1,\ldots,A_{\gamma_{1,m,n}}\}$ be a set of
 $\gamma_{1,m,n}$ distinct matrices.  
 We show that for a generic $L$ $A_1,\ldots, A_d$ are linearly independent.
 Otherwise, for any subspace $L$ of dimension $d$ any $d$ rank one matrices in $L$ must be linearly
 dependent.  (This follows from the fact that linear dependence of $d$ matrices can be stated in terms
 of polynomial equations in the entries of $A_1,\ldots,A_d$.)
 To show that the last condition does not always hold, choose $d$ linearly independent rank one matrices,
 and let $L$ be the subspace spanned by these matrices. 
 
 Assume now that $L$ is a generic subspace of dimension $d\in [(m-1)(n-1)+2,mn-1]$.
 Then $L\cap U_{k,m,n}(\F)$ is a variety of dimension $d-(m-1)(n-1)-1$.
 Similar arguments show that any $d$ generic matrices in $L\cap U_{k,m,n}(\F)$
 are linearly independent.  \qed

 \begin{corol}\label{corexvalgr}
 $\;$
 \begin{enumerate}
 \item\label{corexvalgr1}
 (\ref{baseq}) holds.
 \item
 $\grank(m_1,m_2,(m_1-1)(m_2-1))=(m_1-1)(m_2-1)+1$ for
 $m_1,m_2\ge 2$, i.e. (\ref{(m,n,(m-1)(n-1))}) holds.
 \end{enumerate}
 \end{corol}
 \proof  In view of Proposition \ref{rankeqbm} we discuss first
 the case $m_3 \in[(m_1-1)(m_2-1)+1, m_1m_2-1]$.
 View a generic $\cT=[t_{ijk}]\in\C^{m_1\times m_2\times m_3}$
 as $m_3$ generic matrices $A_k=[t_{ijk}]_{i=j=1}^{m_1,m_2}\in
 \C^{m_1\times m_2}$ for $k=1,\ldots,m_3$.  Hence
 $L=\span(A_1,\ldots,A_{m_3})$ is a generic subspace of
 dimension $m_3$.  Theorem \ref{subspr1m} yields that $L$ is
 spanned by rank one matrices.  Theorem \ref{chartrank} yields
 that $\grank(m_1,m_2,m_3)=m_3$.

 Assume now that $m_3=(m_1-1)(m_2-1)$ and $\cT=[t_{ijk}]
 \in \C^{m_1\times m_2\times m_3}$ be
 a generic tensor.  Let $L\subset \C^{m_1\times m_2}$ be the
 generic subspace defined as above.  Theorem \ref{rnkvarmat} yields that
 $L$ is not spanned by rank one matrices.  Hence the minimal
 dimension of a subspace spanned by rank one matrices
 containing $L$ is at least $m_3+1$.  Let $X\in \C^{m_1\times
 m_2}$ be a generic matrix.  Then $L_1=\span(L,X)$ is a generic
 subspace of dimension $(m_1-1)(m_2-1)+1$.  Hence $L_1$ is
 spanned by rank one matrices.  Therefore $\rank \cT=m_3+1$.
 \qed

 \begin{corol}\label{grank2mn}
 $\grank(2,m_2,m_3)=\min(m_3,2m_2)$  for $2\le m_2\le m_3$.
 \end{corol}

 We now show how to apply the above results to obtain upper estimates of
 $\grank(m_1,m_2,m_3)$ and $\mrank(m_1,m_2,m_3)$.
 Let us start with the case $m_2=m_3\ge 3$.
 \begin{theo}\label{casenmm}  Let $m,n\ge 3$ be integers.
 Then
 \begin{eqnarray}\label{grnmm1}
 \grank(n,m,m)\le \lfloor\frac{n}{2}\rfloor m +
 (n-2\lfloor\frac{n}{2}\rfloor)(m-\lfloor \sqrt{n-1}\rfloor)
 \textrm{ if }
 m \ge 2\lfloor \sqrt{n-1}\rfloor\\
 \grank(n,m,m)\le n(m-\lfloor \sqrt{n-1}\rfloor) \textrm{ if }
 m < 2\lfloor \sqrt{n-1}\rfloor <2(m-1), \label{grnmm2}\\
 \grank(n,m,m)=\min(n,m^2) \textrm{ if } n\ge (m-1)^2+1,
 \label{grnmm3}\\
 \mrank(n,m,m)\le \nonumber\\
 \sum_{i=1}^{\lfloor \sqrt{n-1}\rfloor} (2i-1)(m-i+1)+
 (m- \lfloor \sqrt{n-1}\rfloor^2)(m-\lfloor \sqrt{n-1}\rfloor).
 \label{mrmmm}
 \end{eqnarray}

 \end{theo}

 \proof  We first discuss the $\grank(n,m,m)$.
 Clearly, (\ref{grnmm3}) is implied by Corollary \ref{corexvalgr}.

 Assume now that $n<(m-1)^2+1$, i.e. $2\lfloor \sqrt{n-1}\rfloor
 <2(m-1)$.  Let
 $\tau\in\C^{n\times m\times m}$ be a tensor of the form
 (\ref{tenrep}).  Assume that $(T_{1,1}=[t_{1jk}],\ldots,T_{n,1}
 =[t_{njk}])\in
 (\C^{m\times m})^n$ is a generic point.
 Let $l=\lfloor
 \sqrt{n-1}\rfloor$.  So $n\ge l^2+1$.  Theorem \ref{rnkvarmat} yields that
 $\span (T_{1,1},\ldots,T_{n,1})$ contains at least $\gamma_{m-l,m,m}$
 distinct matrices of rank $m-l$.  It is straightforward to show
 that $\gamma_{m-l,m,m}\ge n$.  Since $(T_{1,1},\ldots,T_{n,1})$ was
 a generic point we may assume $\span (T_{1,1},\ldots,T_{n,1})$
 contain $n$ linearly independent rank $m-l$ matrices $Q_1,\ldots,Q_n$.
 (See the proof of Theorem \ref{subspr1m}.)
 This gives the inequality (\ref{grnmm2}) for
 all $n<(m-1)^2+1$.

 Since $T_{1,1},\ldots,T_{n,1}$ are generic, we can assume that
 $T_{2i-1,1}$ is invertible
 and $T_{2i-1,1}^{-1}T_{2i,1}$ is diagonable.  Hence $T_{2i-1,1},T_{2i,1}$ are contained
 in a subspace spanned by $m$ rank one matrices.
 If $n$ is even we obtain that $\span(T_{1,1},\ldots,T_{n,1})$
 are contained in $\frac{n}{2}m$ dimensional subspace spanned by rank one
 matrices.  Theorem \ref{chartrank}
 yields the inequality (\ref{grnmm1}).  If $n$ is odd, we can
 assume that
 $Q_1=T_{n,1}-\sum_{i=1}^{\lfloor\frac{n}{2}\rfloor} T_{i,1}$
 has at rank $m-\lfloor \sqrt{n-1}\rfloor$.  Hence, we deduce
 (\ref{grnmm1}) in this case too.

 We now prove the inequality (\ref{mrmmm}).
 We assume the worst case which will give the upper bound.
 So it is enough to consider the case where  $T_{1,1},T_{2,1},\ldots,
 T_{n,1}$ linearly independent.  Now we choose a  new base
 $S_1,\ldots,S_n$ in  $\span (T_{1,1},\ldots,T_{n,1})$
 such that $\rank S_1\ge \rank S_2\ge \ldots\ge \rank S_n$.
 So the worst case is $\rank S_1=m$.  Since any $2$ dimensional
 space contains a singular matrix we can assume that $\rank S_i\le
 m-1$ for $i=2,3,4$.  According to Theorem \ref{rnkvarmat}
 any $5$ dimensional vector space contains a nonzero matrix of
 rank $m-2$ at most.  Hence $\rank S_i\le m-2$ for $i=5,6,7,8,9$.
 Theorem \ref{rnkvarmat} implies that any subspace of dimension
 $10$ contains a nonzero matrix of rank $m-3$.  Hence $\rank S_i \le
 m-3$ for $i=10,\ldots,$.  Continuing the use of Theorem
 \ref{rnkvarmat},
 and combing it with Theorem \ref{chartrank} we deduce (\ref{mrmmm}).
 \qed

 Use Corollary \ref{grankub}, Proposition \ref{incrgrank} and the above theorem to deduce:
 \begin{corol}\label{basinmmm}
 \begin{eqnarray*}
 &4 \le\grank(3,3,3)\le 5=1\cdot 3+2, & \mrank(3,3,3)\le 7=3+2+2,\\
 & \grank(4,3,3)=5 \;(4=(3-1)^2), & \mrank(4,3,3)\le 9=3+2+2+2,\\
 &\grank(5,3,3)=5 \;(5> (3-1)^2), & \mrank(5,3,3)\le 10=3+2+2+2+1,\\
 &6 \le \grank(3,4,4)\le 7=1\cdot 4+3, & \mrank(3,4,4)\le 10=4+3+3,\\
 & 7 \le \grank(4,4,4)\le 8=2\cdot 4, & \mrank(4,4,4)\le 13=4+3+3+3,\\
 & 8 \le \grank(5,4,4)\le 10=2\cdot 4+2, & \mrank(5,4,4)\le 15=4+3+3+3+2,\\
 & 7\le \grank(3,5,5)\le 9=1\cdot 5+4, &\mrank (3,5,5)\le 13=5+4+4,\\
 & 9\le \grank(4,5,5)\le 10=2\cdot 5, &\mrank (4,5,5)\le 17=5+4+4+4,\\
 & 10\le \grank(5,5,5)\le 13=2\cdot 5+3, &\mrank (5,5,5)\le 20=5+4+4+4+3.
 \end{eqnarray*}
 \end{corol}

 Recall that in all the examples of $\grank(n,m,m)$ given by Corollary \ref{basinmmm}  we know that $\grank(3,3,3)=5, \grank(3,5,5,)=8$,
 while all other values of $\grank(n,m,m)$ are given by the lower bound.
 It is claimed that $\mrank(3,3,3)=5$ \cite{Roc93}.

 Note that if $n$ is even and $m\gg n$ then the upper bound (\ref{grnmm1}) combined with
 Corollary \ref{grankub} implies that $\grank(n,m,m)$ is of order $\frac{nm}{2}$.
 However if $n=O(m^{1+a})$ for $a\in (0,1]$  then the upper bounds (\ref{grnmm1}--\ref{grnmm2})
 are not of the right order, (which is $m^2$).

 \section{Typical ranks of real $3$-tensors}

 The study of the rank of a real $3$-tensor is closely related to the real
 semi-algebraic geometry.  See \S\ref{Sec:AppendixR} for the results in semi-algebraic geometry needed
 here.
 \begin{theo}\label{grnankreal}
 The space $\R^{m_1\times m_2\times m_3}, m_1,m_2,m_3\in\N$, contains a
 finite number of open connected disjoint semi-algebraic sets
 $O_1,\ldots, O_M$ satisfying the following properties.
 \begin{enumerate}\label{grnankreal1}
 \item $\R^{m_1\times m_2\times m_3}\backslash \cup_{i=1}^M O_i$
 is a closed semi-algebraic set $\R^{m_1\times m_2\times
 m_3}$ of dimension strictly less than $m_1m_2m_3$.
 \item\label{grnankreal2} Each $\cT\in O_i$ has rank $r_i$ for
 $i=1,\ldots,M$.
 \item\label{grnankreal3}
 $\min (r_1,\ldots, r_M)=\grank(m_1,m_2,m_3)$.
 \item\label{grnankreal4}
 $\rmrank(m_1,m_2,m_3):=\max (r_1,\ldots, r_M)$
 is the minimal $k\in\N$ such that the closure of
 $\f_k((\R^{m_1}\times \R^{m_2}\times
 \R^{m_3})^k)$ is equal to $\R^{m_1\times m_2\times m_3}$.
 \item\label{grnankreal5}  For each integer $r\in [\grank(m_1,m_2,m_3),\rmrank(m_1,m_2,m_3)]$
 there exists $r_i=r$ for some integer $i\in [1,M]$.
 \end{enumerate}
 \end{theo}
 \proof Consider the
 polynomial map $\f_k:(\C^{m_1}\times \C^{m_2}\times \C^{m_3})^k\to
 \C^{m_1\times m_2\times m_3}$ be given by (\ref{deffk}).
 Note that $\f_k:(\R^{m_1}\times \R^{m_2}\times \R^{m_3})^k\to
 \R^{m_1\times m_2\times m_3}$.  Denote by $Y_k$ and $Q_k$
 the closure of $\f_k((\C^{m_1}\times \C^{m_2}\times
 \C^{m_3})^k)$ and  $Z_k:=\f_k((\R^{m_1}\times \R^{m_2}\times
 \R^{m_3})^k)$ respectively.  Clearly,
 $$Y_i\subseteq Y_{i+1},\; Q_i\subseteq Q_{i+1} \textrm{ for
 }i\in\N, \; Y_{m_1m_2m_3}=\C^{m_1\times m_2\times m_3},\;
 Q_{m_1m_2m_3}=\R^{m_1\times m_2\times m_3}.$$
 Let $\rmrank(m_1,m_2,m_3)$ be the smallest $k$ such that
 $Q_k=\R^{m_1\times m_2\times m_3}$.

 Let $q=\grank(m_1,m_2,m_3)$.  Then $Y_{q-1}$ is a strict
 complex subvariety of $\C^{m_1\times m_2\times
 m_3}$.  (See Definition \ref{defrkXk}.)  In particular $Y_{q-1}^{\R}=
 Y_{q-1}\cap \R^{m_1\times m_2\times m_3}$ is a strict real subvariety
 of $\R^{m_1\times m_2\times m_3}$.  Hence $Q_{q-1}\subseteq
 Y_{q-1}^{\R}$ is a semi-algebraic of dimension $\dim Y_{q-1}$ at most,
 which is strictly less than $m_1m_2m_3-1$.  In particular
 \begin{equation}
 \rmrank(m_1,m_2,m_3)\ge
 \grank(m_1,m_2,m_3).
 \label{inrmrnkgrnk}
 \end{equation}
 From the proof of Theorem \ref{maxrank} it follows that
 there exists an algebraic subset $X_q \subset (\C^{m_1}\times \C^{m_2}\times
 \C^{m_3})^q$ such that $\rank \rD\f_q$ is $m_1m_2m_3$ at each
 point of $(\C^{m_1}\times \C^{m_2}\times \C^{m_3})^q\backslash
 X_q$.  Then $X^\R_{q}=X_q\cap (\R^{m_1}\times \R^{m_2}\times
 \R^{m_3})^q$ is a real algebraic set of $(\R^{m_1}\times \R^{m_2}\times
 \R^{m_3})^q$.  Thus the Jacobian of the real map $\f_q:
 (\R^{m_1}\times \R^{m_2}\times \R^{m_3})^k \to \R^{m_1\times
 m_2\times m_3}$ has rank $m_1m_2m_3$ at each point of the
 open semi-algebraic set
 $P_q:=(\R^{m_1}\times \R^{m_2}\times \R^{m_3})^q\backslash
 X^\R_q$.  Hence $\f_q(P_q)$ is an open semi-algebraic set in
 $\R^{m_1\times m_2\times m_3}$.  Therefore
 $\f_q(P_q)\backslash Y^\R_{q-1}$ is an open semi-algebraic set in
 $\R^{m_1\times m_2\times m_3}$. Clearly
 $Q_q\setminus Q_{q-1}\supseteq Q_q\setminus Y^\R_{q-1}\supseteq \f_q(P_q)\backslash Y^\R_{q-1} $.
 Hence the interior of $Q_q\setminus Q_{q-1}$, denoted as $\inter (Q_q\setminus Q_{q-1})$ is an open semi-algebraic set, which consists of tensors
 of rank $q$ exactly.
 The theory of semi-algebraic sets implies that $\inter (Q_q\setminus Q_{q-1})=\cup_{i=1}^{M_1} O_i$, where each $O_i$ is an open semi-algebraic set.  Observe next that the semi-algebraic set
 $(Q_q\setminus Q_{q-1})\setminus \inter (Q_q\setminus Q_{q-1}$) has dimension $m_1m_2m_3-1$ at most.
 Since $\dim Q_{q-1}\le m_1m_2m_3-1$ we deduce that
 \begin{equation}\label{QqminClosin}
 \dim Q_q\setminus \textrm{Closure} (\cup_{i=1}^{M_1} O_i)\le m_1m_2m_3-1.
 \end{equation}
 Suppose $Q_q=\R^{m_1\times m_2\times m_3}$, i.e. equality holds in
 (\ref{inrmrnkgrnk}), so $M=M_1$.  We claim that $W_q:=\R^{m_1\times m_2\times m_3}
 \setminus \textrm{Closure} (\cup_{i=1}^{M_1} O_i)$ is an empty set. Otherwise $W_q$
 is a nonempty open semi-algebraic set.  Hence $\dim W_q=m_1m_2m_3$ which contradicts
 (\ref{QqminClosin}).  The proof of the theorem is completed in this case.

 Assume now that $Q_q \subsetneq \R^{m_1\times m_2\times
 m_3}$.  Recall that $\dim \textrm{Closure} (S)\setminus S < \dim S$ for any semi-algebraic set.
 Hence $\dim Q_{q+1}=\dim Z_{q+1}$.  We claim that
 $\dim (Z_{q+1}\setminus Q_q)=m_1m_2m_3$, i.e. the interior of
 $Z_{q+1}\setminus Q_q$ contains an open set.  Assume to the contrary that
 that $\dim (Z_{q+1}\setminus Q_q)<m_1m_2m_3$.  Hence $\dim (Z_{q+1}\setminus Z_q)<m_1m_2m_3$,
 ($\dim Q_q\setminus Z_q<\dim Z_q=m_1m_2m_3$.)  So a sum of generic $q+1$ real rank one tensors
 is a sum of generic $q$ real rank one tensors.  Hence a sum of generic $m_1m_3m_3$ rank one tensors
 is a sum of $q$ generic rank one tensors.  So $Q_q=\R^{m_1\times m_2\times m_3}$, which contradicts our assumption.
 Thus, the interior of $Q_{q+1}\setminus Q_q$ is an open semi-algebraic set, which is a union of
 disjoint open connected semi-algebraic sets $O_{M_1+1},\ldots,O_{M_2}$.  Note that the rank $\cT\in O_j$ is $\grank(m_1,m_2,m_3)+1$ for $j=M_1+1,\ldots,M_2$.  Continue in this manner we deduce
 the rest of the theorem.
 \qed
 \begin{defn}\label{deftypbrnk}
 Let $r$ be a positive integer. $\cT\in \R^{m_1\times m_2\times m_3}$ has a border rank $r$,
 denoted as $\brank \cT$,
 if $\cT\in \mathrm{Closure}\; \f_r((\R^{m_1}\times \R^{m_2}\times \R^{m_3})^r)\setminus
 \mathrm{Closure}\; \f_{r-1}((\R^{m_1}\times \R^{m_2}\times \R^{m_3})^{r-1})$.
 ($\f_0((\R^{m_1}\times \R^{m_2}\times \R^{m_3})^0=\{0\}$.)
 $r$ is called an $(m_1,m_2,m_3)$ \emph{typical} rank, or simply typical rank, if
 $r\in [\grank(m_1,m_2,m_3),\rmrank(m_1,m_2,m_3)]$.
 \end{defn}

 The proof of Theorem \ref{grnankreal} yields.
 \begin{corol}\label{probtypbord}  Assume that the entries of $\cT\in \R^{m_1\times m_2\times m_3}$
 are independent random variables with standard normal Gaussian distribution.
 Then the probability that $\rank \cT=r$ is positive if and and only if $r$ is a typical rank.
 Assume that $r$ is a typical rank.  Then the probability that $\rank \cT> \brank\cT$, provided that $(\rank \cT -r)(\brank \cT-r)=0$,  is $0$.
 In particular, the probability that $\rank\cT=\grank(m_1,m_2,m_3)$ is positive.
 \end{corol}

 The last part of this Corollary is shown in \cite[Appendix B]{tBS06} for $m_1=m_2=4, m_3=3$.
 For $l=2\le m\le n$ the following is known:
 $\rmrank(2,m,m)=\grank(2,m,m)+1=m+1$ \cite{tB91} and $\rmrank(2,m,n)=\grank(2,m,n)=\min(n,2m)$
 for $m<n$ \cite{tBK99}.
 \cite{Roc93} claims that $\rmrank(3,3,3)=\grank(3,3,3)=5$. It is shown in \cite{tB04} that
 $\rmrank(3,3,5)=\grank(3,3,5)+1=6$.
 For other additional known results for typical rank see \cite{CBLC}.  In particular,
 $\rmrank(4,4,12)=\grank(4,4,12)+1=12$ \cite[Table I]{CBLC}.
 We now give additional examples, where a strict inequality holds in
 (\ref{inrmrnkgrnk}).  All of them, except the above mentioned examples, are new.
 \begin{theo}\label{examinermaxgr}
 In the following cases
 $\rmrank(m_1,m_2,m_3)>\grank(m_1,m_2,m_3)$.
 \begin{enumerate}
 \item \label{examinermaxgr1} $m_1=m_2=m\ge 2, m_3=(m-1)^2+1$.
 \item \label{examinermaxgr2} $m_1=m_2=4, m_3=11,12$.
 \end{enumerate}
 \end{theo}

 We do not know if
 $\rmrank(m,m,(m-1)^2+1)=\grank(m,m,(m-1)^2+1)+1$ for $m\ge 4$.
 To prove Theorem \ref{examinermaxgr} we need a few auxiliary results.
 The following result is known, e.g. \cite[Proposition 5.2]{FK}.
 \begin{prop}\label{irw2pc}  Let $\F=\C,\R$, $n\ge 2, p\ge 1$ be integers and
 assume that $p\le \lfloor \frac{n}{2}\rfloor$.
 Let $\proj \rA_n(\F)\supseteq\proj \rW_{2p,n}(\F)$ be the projective variety of all
 (nonzero) skew symmetric matrices and the projective subvariety of all skew
 symmetric matrices of rank $2p$ at most respectively.  Then $\proj
 \rW_{2p,n}(\F)$
 is an irreducible projective variety in
 $\proj \rA_n(\F)$ of codimension $n-2p \choose 2$.
 The variety of its singular points is $\proj \rW_{2(p-1),n}(\F)$.
 \end{prop}
 \begin{corol}\label{nonitsucsks}  A generic subspace $L$ of
 the linear space of $n\times n$ skew symmetric matrices
 $\rA_n(\F)\subset \F^{n\times n}$ of dimension ${n-2p \choose 2}$
 does not contain a nonzero matrix of rank $2p$ at most.
 In particular, for each generic point
 $\T:=(T_1,\ldots,T_{n-2p \choose 2})\in \rA_n(\F)^{n-2p \choose
 2}$, there exists an open neighborhood of $O\subset
 \rA_n(\F)^{n-2p \choose 2}$
 such that for each $\X:=(X_1,\ldots,X_{n-2p \choose 2})\in O$,
 $L(\X):=\span(X_1,\ldots,X_{n-2p \choose 2})$ is a subspace of
 dimension of ${n-2p \choose 2}$ which does not contain a
 nonzero matrix of rank $2p$ at most.
 \end{corol}
 \proof  A subspace $L\subset \rA_n(\F)$ of dimension $d$ induces
 a linear space $\proj L$ of dimension $d-1$ in the projective
 space $\proj \rA_n(\F)$. Hence the dimension count implies
 that  $\proj L\cap \P \rW_{2p,n}(\F)=\emptyset$ for a generic subspace
 $L$ of dimension ${n-2p \choose 2}$.  Hence $L$ does not
 contain a nonzero matrix of rank $2p$ at most.

 A generic point $\T\in \rA_n(\F)^{n-2p \choose 2}$ generates a
 generic subspace $L(\T)$ of dimension ${n-2p \choose 2}$.
 Hence $\proj L(\T)\cap \proj\rW_{2p,n}(\F)=\emptyset$.  For a small
 enough open neighborhood $O$ of $\T$, for any $\X\in O$,
 the subspace $L(\X)$ is a perturbation of $L(\T)$.
 Hence $\proj L(\X)\cap \proj\rW_{2p,n}(\F)=\emptyset$.  \qed

 It is well known that for $\F=\R$ the above corollary can be
 improved for certain values of $n,p$.  See \cite{FK} and the
 references therein.  We now bring a well known improvement
 of the above corollary for $n=4, p=1$.

 \begin{prop}\label{case42}  There exists an neighborhood
 $O$ of $\T=(T_1,\ldots,T_l)\in \rA_4(\R)^l$ such that
 for any $\X=(X_1,\ldots,X_l)\in \rA_4(\R)^l$ the subspace
 $L(\X)$ does not contain a matrix of rank $2$ for $l=2,3$.
 \end{prop}
 \proof  Let $l=3$ and
 $$T_1=\left[ \begin{matrix}
 0&1&0&0\\-1&0&0&0\\0&0&0&1\\0&0&-1&0\end{matrix}\right],
 T_2=\left[ \begin{matrix}
 0&0&1&0\\0&0&0&-1\\-1&0&0&0\\0&1&-0&0\end{matrix}\right],
 T_3=\left[ \begin{matrix}
 0&0&0&1\\0&0&1&0\\0&-1&0&0\\-1&0&0&0\end{matrix}\right].
 $$
 Let $\T=(T_1,T_2,T_3)$.  Note that any nonzero matrix $B\in
 L(\T)$ is a multiple of an orthogonal matrix.  Hence $\rank
 B=4$ and $\dim L=3$.  Thus $\proj L(\T)\cap \proj \rW_{2,4}(\R)=\emptyset$.
 Therefore, there exists a small open neighborhood $O$ of
 $\T$ such that for any $\X=(X_1,X_2,X_3)\in O$ $\proj L(\X)\cap \proj \rW_{2,4}(\R)=\emptyset$.

 Similar results hold for $l=2$ if we let $\T=(T_1,T_2)$.  \qed
 The next result appears in \cite{FL}.
 \begin{prop}\label{frloew}
 Let $\rS_{n,0}\subset \R^{n\times n}$ be the subspace of real
 symmetric matrices of trace zero.  Then $\rS_{n,0}$ is an
 $\frac{(n+1)n}{2}-1$ dimensional subspace which does not
 contain a rank one matrix.
 \end{prop}
 \proof  Clearly, $\dim \rS_{n,0}=\frac{(n+1)n}{2}-1$.  Assume to the
 contrary that a rank one matrix $B$ is in $\rS_{n,0}$.  Since $B$ is
 symmetric $B=\pm\x\x\trans$, where $\0\ne\x\in\R^n$.  Then trace
 $B=\pm\x\trans \x=0$.  So $\x=\0$, contradicting our assumption.
 \qed \\
 \textbf{Proof of Theorem \ref{examinermaxgr}}.
 We first begin with the case $(m,m,l=(m-1)^2+1)$.
 Assume first $m=2,3$.  Note that $\dim \rS_{m}=l$.
 Choose a basis $T_1,\ldots,T_l$ in $\rS_m$.
 Let $\T=(T_1,\ldots,T_l)$. Proposition \ref{frloew}
 yields that $\proj\rS_{n,0}\cap \proj U_{1,m,m}=\emptyset$.
 The arguments of the proof of Corollary \ref{nonitsucsks}
 yield that there exists an open neighborhood $O$ of $\T\in
 (\R^{m\times m})^l$ so that for each $\X=(X_1,\ldots,X_l)\in
 (\R^{m\times m})^l$ we have $\proj L(\X)\cap \proj U_{1,m,m}=\emptyset$.
 Hence $L(\X)$ is not spanned by rank one matrices.

 Let $\cT=[t_{ijk}]\in \R^{m\times m\times l}$ be the set of
 $C\subset \R^{m\times m\times l}$ of
 all $3$-tensors such that $\X\in O$, where
 $X_k:=[t_{ijk}]_{i=j=1}^m$ for $k=1,\ldots,l$.
 Clearly, $C$ is open.
 Theorem \ref{chartrank} implies that the $\mathrm{rank}_{\R} \cT>l$ for each
 $\cT\in C$.  In view of  Theorem \ref{grnankreal},
 $C$ has a nontrivial intersection with at least
 one $O_i$.  Hence $r_i>l=\grank(m,m,l)$.

 Assume now that $m>3$.  Let $L_1\subset \rA_m(\R)$ be a generic
 subspace of dimension $m-2\choose 2$.  Then $L_1$ does not
 contain a matrix of rank $2$. Clearly $\rS_{m,0}\cap
 L_1=\{0_{m\times m}\}$.  Then $L=S_{m,0}+L_1$ is $l=(m-1)^2+1$
 dimensional subspace of trace zero matrices.
 Observe that if $B\in L$ then $B\trans\in L$.   We claim that $L$ does not contain a
 rank one matrix $B\in \R^{m\times m}$.  Assume to the contrary
 that $B\in L$ is a rank one matrix.   Proposition \ref{frloew}
 implies that $B\not\in \rS_{m,0}$.  So
 $$B=B_1+B_2,\;
 B_1=\frac{1}{2}(B+B\trans)\in\rS_{m,0},
 B_2=\frac{1}{2}(B-B\trans)\in L_1.$$
 Since $B$ is a rank one nonsymmetric matrix $B_2$ is a skew
 symmetric matrix of rank $2$.  This contradicts our
 assumption.  Hence $\proj L\cap \proj U_{1,m,m}=\emptyset$.
 The above arguments show that
 $\rmrank(m,m,l)>l=\grank(m,m,l)$.

 Assume finally that $m=4$ and $l=11,12$.
 Repeat the above arguments where $L_1$ has dimension $2$ or
 $3$, as given in Proposition \ref{case42}.  \qed

 \section{Appendix: Complex and real algebraic geometry}\label{Sec:Appendix}
 In this section we give basic facts in complex and real algebraic geometry needed for this paper.
 The emphasize is on simplicity and intuitive understanding.  We supply references for completeness.
 Our basic references are \cite{Mum76}, \cite{Sha77} and \cite{Har} for complex algebraic geometry, and \cite{BCR} for
 real algebraic geometry.

 We first start with some general definitions which hold for general field $\F$.
 Denote by $\F[x_1,\ldots,x_n],\F(x_1,\ldots,x_n)$ the ring of polynomials and its field
 of rational functions in $n$ variables $x_1,\ldots,x_n$ with coefficients in $\F$ respectively.  We will identify $\F[\x]=\F[x_1,\ldots,x_n], \F(\x)=\F(x_1,\ldots,x_n)$, where $\x=(x_1,\ldots,x_n)\trans\in\F^n$.
 For $p_1,\ldots,p_m\in\F[\x]$ denote by $Z(p_1,\ldots,p_m)=\{\y\in\F^n,\;p_i(\y)=0, i=1,\ldots,m\}$.
 Equivalently let $\mP=(p_1,\ldots,p_m)\trans$ be a polynomial map $\mP:\F^n \to \F^m$.  Then $Z(p_1,\ldots,
 p_m)=\mP^{-1}(\0)$.  $V\subset\F^n$ is called an \emph{algebraic set}, if $V=Z(p_1,\ldots,p_m)$ for some $p_1,\ldots,p_m\in\F[\x]$.  Note that $\emptyset$ and $\F^n$
 algebraic sets.

 Recall that $\bP\F^n$, the $n$-dimensional projective space over $\F$, is identified with one dimensional subspaces of $\F^{n+1}$, i.e. lines through the origin in $\F^{n+1}$.
 So $\F^n$ is viewed as a subset of $\bP\F^n$ where each $\x=(x_1,\ldots,x_n)\trans$ is identified with
 a one dimensional subspace spanned by $\hat x=(x_1,\ldots,x_n,1)\trans$.  $\bP\F^n$ can be viewed as the
 union of two disjoint sets $\F^n$ and $\bP\F^{n-1}$, where $\bP\F^{n-1}$ is all one dimensional subspaces
 in $\F^{n+1}$ spanned by nonzero $\y=(y_1,\ldots,y_n,0)\trans$.

 Denote by $\F_h[\y], \y=(y_1,\ldots,y_{n+1})\trans$, the set of homogeneous polynomials in $y_1,\ldots,y_{n+1}$.  Let $q_1,\ldots,q_m\in \F_h[\y]$.  Consider the variety $Z(q_1,\ldots,q_m)\subset \F^{n+1}$.   If $\0\ne \y\in Z(q_1,\ldots,q_m)$ then $\span(\y)\subset Z(q_1,\ldots,q_m)$.
 Hence $Z(q_1,\ldots,q_m)$ induces a subset $\tilde Z(q_1,\ldots,q_m)\subset \bP\F^n$.  (If
 $Z(q_1,\ldots,q_m)=\{\0\}$ then $\tilde Z(q_1,\ldots,q_m)=\emptyset$.)  $V\subseteq \bP\F^{n+1}$
 is called a \emph{projective algebraic set} if $V=\tilde Z(q_1,\ldots,q_m)$ for some $q_1,\ldots,q_m\in\F_h[\y]$.
 It is easy to show that an intersection and union of two affine or projective algebraic sets is an affine or projective algebraic.  An affine or projective algebraic set is called irreducible if it cannot be written as the union of two proper algebraic subsets.  An irreducible affine or projective algebraic set is called an affine or projective variety respectively.  (An affine variety will be referred sometimes as variety.)
 Let $V$ be a projective variety in $\bP\F^n$, and $W\subsetneq V$ a projective algebraic set.  Then $V\setminus W$ is called a \emph{quasi-projective} variety.   Note that an affine variety $Z(p_1,\ldots,p_m)$ can be viewed as a quasi projective variety.  First homogenize $p_1,\ldots,p_m$ to $\hat p_1,\ldots,\hat p_m\in\F_n[\y]$.  Let $W\subset \bP\F^n$ to be the zero set of $y_{n+1}=0$.  Then
 $Z(p_1,\ldots,p_m)$ can be identified with $Z(\hat p_1,\ldots,\hat p_m)\backslash W$.

 \subsection{Complex algebraic sets and polynomial maps}\label{Sec:AppendixC}
 In this section $\F=\C$.
 Let $\mP=(p_1,\ldots,p_m):\C^n\to\C^m$ be a polynomial map.
 Denote by $\rD\mP(x)$, the derivative of $\mP$ or the Jacobian matrix of $\mP$,
 the matrix $[\frac{\partial p_i}{\partial x_j}]_{i=j=1}^{m,n}$.
 For any $U\subseteq \C^n$
 denote $\mathrm{rank}_U \rD\mP=\max_{\x\in U} \rank \rD\mP(\x)$.
 Assume that $U$ is a variety.
 Note that the set $\Sing U=\{\x\in U,\;\rank \rD\mP(\x)<\mathrm{rank}_U \rD\mP\}$ is a strict algebraic subset of $U$.
 (Observe that $\x\in\Sing U$ if and only if all minors of $ \rD\mP(\x), \x\in U$ of order $\mathrm{rank}_U \rD\mP$ vanish.)
 $\Sing U$ is called the set of singular points of $U$.  Let $V=Z(p_1,\ldots,p_m)$ be a variety.
 The dimension of $V$, denoted by $\dim V$, equals to $n-\mathrm{rank}_V \rD\mP$.  Then $V\backslash \Sing V$, the set of regular (smooth) points of $V$, is a quasi-projective variety, and a complex manifold of dimension $\dim V$.  See \cite[\S1A]{Mum76}.  For any variety $V$ and a strict algebraic subset $W$ in $V$, the quasi-projective variety $V\setminus W$ is connected \cite[Cor 4.16]{Mum76}, and its dimension equal to the dimension of the complex manifold $V\setminus (W\cup \Sing V)$, which is $\dim V$.
 We say that a given property holds \emph{generically} in $V$, if it holds for each $\x\in V\setminus W$,
 for some strict algebraic subset $W$ of $V$, where $W$ depends on the given property.

 Hilbert basis theorem, (Nullstellensatz),
 claims that a countable intersection of algebraic sets is an algebraic set \cite[p'17]{Sha77}.
 An algebraic set $U\subset \C^n$ is a union of finitely many pairwise distinct varieties $U_1,\ldots,U_k$, and this decomposition is unique \cite[Thms I.3.1, I.3.2]{Sha77}.  We define $\dim U=\max \dim U_i$.
 A product of two irreducible varieties is an irreducible variety \cite[Thm I.3.3]{Sha77}.
 Similar results holds for projective algebraic sets.

 A set $V\subset \C^{n}$ is called a \emph{constructible} algebraic set of dimension $d$ if it can be represented as $V\setminus W$ were $V$ is an algebraic set of dimension $d$ and $W$ is a constructible algebraic set of dimension $d-1$ at most \cite{Har}.  Note that a constructible algebraic set of dimension $0$ is a set consisting of a finite number of points.  It is easy to show that a finite union and a finite intersection of constructible algebraic sets is a constructible algebraic set.  Finally if $V,W\subset \C^n$ are constructible algebraic sets then $V\setminus W$
 is constructible algebraic.

 Let $\mP$ be a polynomial map as above.  From the definition of an algebraic set we deduce that for any
 algebraic set $W\subset \C^m$ the set $\mP^{-1}(W)$ is an algebraic set of $\C^n$.
 Denote $\rank \rD\mP=\mathrm{rank}_ {\C^n}\rD\mP$.  Then
 $V=\textrm{Closure }\mP(\C^n)$ is a variety, of dimension $\rank\rD\mP$.
 (Here the closure is in the standard topology in $\C^n$ or $\R^n$.)
 Moreover, $\Sing \mP= \{\x\in\C^n,\; \rank \rD\mP(x)<\rank\rD\mP\}$ is a strict algebraic subset of
 $\C^n$.
 Hence $\mP(\C^n\setminus \Sing \mP)$ is a constructive algebraic variety in $\C^m$ of dimension
 $\rank \rD\mP$ \cite{Har}. Furthermore, there exists a strict algebraic set $W\subsetneq V$, such that for each
 $\z\in V\setminus W$ the algebraic set $\mP^{-1}(\z)$ is a disjoint union of $k$ varieties $U_1(z),\ldots,
 U_k(z)\subset \C^n$, each of dimension $n-\rank \rD\mP$.  The integer $k$ is independent of $\z\in V\setminus W$, and is called the \emph{degree} of $\mP$ \cite[Corol. 3.15-3.16]{Mum76}.

 More general, let $U\subset \C^n$ be a constructible algebraic set.  Then $\mP(U)\subset \C^m$ is a constructible algebraic set of dimension $\mathrm{rank}_U \rD\mP$.  This applies in particular to a projections $\mP$, where $\mP(\x)$ obtained from $\x$ be deleting a number of coordinates.  See \cite[\S3-4]{Sha77}.

 \subsection{Real semi-algebraic sets and polynomial maps}\label{Sec:AppendixR}
 In this subsection the topology on $\R^n$ is assumed to be the standard topology: open sets, closed sets, the interior and the closure of sets are in the standard topology of $\R^n$.
 A real algebraic set in $\R^n$ is the zero set of $m$ polynomials $p_1,\ldots,p_m\in \R[\x]$, and
 is denoted by $Z^{\R}(p_1,\ldots,p_m)\subset\R^n$.
 We can view $p_1,\ldots,p_m$ as polynomials with complex variables $\z=(z_1,\ldots,z_n)\trans\in\C^n$
 with real coefficients.  Then
 $U=Z(p_1,\ldots,p_m)=\{\z\in\C^n,\;p_1(\z)=\ldots=p_m(\z)=0\}$ and $U^{\R}=Z^{\R}(p_1,\ldots,p_m)=U\cap\R^n$.  $Z^{\R}(p_1,\ldots,p_m)$ is called irreducible,
 if $Z(p_1,\ldots,p_m)$ is irreducible.  Since any algebraic set $U\subset \C^n$ is a finite union
 of pairwise distinct irreducible varieties $V_1,\ldots,V_k$ it follows that any real algebraic set
 is a finite union of irreducible real algebraic sets.  A set $S$ is called \emph{semi-algebraic}
 if $S$ is a finite union of subsets $S_1,\ldots,S_k$, where each $S_i$ is of the following form.  There exists an algebraic set $V_i^{\R}\subset \R^n$ and a finite number of polynomials $g_{1,i},\ldots g_{n_i,i}\in\R[\x]$ such that
 $S_i=\{\x\in V_i^{\R},\; g_{j,i}(\x)> 0,\;j=1,\ldots,n_i\}$ for $i=1,\ldots,k$.  Here each $n_i\ge 0$.  So if $n_i=0$ then $S_i=V_i^{\R}$.
 (Algebraic set is semi-algebraic.) Since each algebraic set is a finite union of irreducible real varieties we may assume that in the
 definition of semi-algebraic set $S$ each $V_i^{\R}$ is irreducible.  Furthermore, without loss of generality, we may assume that
 each $S_i\subset V_i^{\R}$ is relative open, i.e. $S_i$ is a nonempty intersection of an open set in $\R^n$ and $V_i^{\R}$.
 Hence $\dim S_i=\dim V_i^{\R}$, and $\dim S=\max \dim S_i$.
 See \cite[\S2.8]{BCR}.

 Semi-algebraic sets are stable under finite union, finite intersection, taking complements and closures  \cite[\S2.2]{BCR}.  (I.e.  all the above operations on semi-algebraic sets yield semi-algebraic sets.)  Hence if $S,T$ are semi-algebraic subsets of $\R^n$ then $A\setminus B=A\cap (\R^n\setminus B)$ is a semi-algebraic set.   For any semi-algebraic set $S$ the following inequality holds $\dim \textrm{Closure} (S)\setminus S < \dim S$ \cite[Prop 2.8.13]{BCR}.

 A projection of semi-algebraic set is
 semi-algebraic \cite[Thm 2.2.1]{BCR}.  Hence the image of a semi-algebraic set by a polynomial map
 is semi-algebraic \cite[Prop 2.2.7]{BCR}.  The closure and the interior of semi-algebraic
 set are semi-algebraic \cite[Prop2.2.2]{BCR}.
 Every open semi-algebraic subset $S$ of $\R^n$ is a finite union of disjoint open connected semi-algebraic sets in $\R^n$.  For more general statement see \cite[Thm 2.4.4]{BCR}.


\begin{thebibliography}{99}

 \bibitem{AOP} H. Abo, G. Ottaviani and C. Peterson, Induction for secant varieties of Segre
 varieties, \emph{Trans. Amer. Math. Soc.}, 361 (2009), 767--792.

 \bibitem{BCLR} D. Bini, M. Capovani, G. Lotti, and F. Romani, $O(n^{2.7799})$ complexity for matrix
 multiplication, \emph{Inf. Proc. Letters} 8 (1979), 234--235.

 \bibitem{BCR} J. Bochnak, M. Coste and M.F. Roy,
 \emph{Real algebraic geometry}, Ergebnisse der Mathematik und ihrer Grenzgebiete (3),
 36, Springer-Verlag, Berlin, 1998.


 \bibitem{BCS} P. B\"urgisser, M. Clausen, M.A. Shokrollahi,
 \emph{Algebraic complexity theory}, with the collaboration of Thomas Lickteig,
 \emph{Grundlehren der Mathematischen Wissenschaften},
 315. Springer-Verlag, Berlin, 1997. xxiv+618 pp.

 \bibitem{CC} J. Carrol and J. Chang, Analysis of individual
 differences in multidimensional scaling via N-way generalization of
 "Eckhart-Young" decomposition, \emph{Psychometrika} 9 (1970),
 267-283.

 \bibitem{CGG} M.V. Catalisano, A.V. Geramita, A. Gimigliano, Ranks of tensors,
 secant varieties of Segre varieties and fat points,
 \emph{Linear Algebra Appl.}  355  (2002), 263--285.

 \bibitem{CBLC} P. Comon, J.M.F. Ten Berge, L. De Lathauwer and J. Castaing,
 Generic and typical ranks of multi-way arrays, \emph{Linear Algebra Appl.}
 430 (2009), 2997-3007.

 \bibitem{CB} R. Coppi and S. Bolasco, Editors, \emph{
 Multiway Data Analysis}, Elsevier Science Publishers, (North-Holland), 1989.

 \bibitem{DL08} V. De Silva and L-H Lim, Tensor rank and the ill-posedness of the best low-rank approximation problem, \emph{SIAM J. Matrix Anal. Appl.} 30 (2008) 1084--1127.

 \bibitem{Fr1} S. Friedland, 3-Tensors: ranks and approximations,
 Workshop on Algorithms for Modern Massive Data Sets,
 Stanford-Yahoo, June 21-24, 2006,
 www.math.uic.edu/$\sim$friedlan/yahostan06.pdf
 \bibitem{FK} S. Friedland and C. Krattenthaler,
 $2$-adic valuations of certain ratios of products of
 factorials and applications, \emph{ Linear Algebra Appl.},
 426 (2007), 159-189, (arXiv.math.NT/0508498 v2)
 \bibitem{FL} S. Friedland and R. Loewy,  Subspaces of symmetric matrices
 containing matrices with a multiple first eigenvalue, \emph{Pacific J. Math.}
 62 (1976), 389-399.

 \bibitem{FMMN} S. Friedland, V. Mehrmann, A. Miedlar and
 M. Nkengla, Fast low rank approximations of matrices and
 tensors, 2008, www.matheon.de/preprints/4903.

 \bibitem{Gan} F.R. Gantmacher, \emph{The Theory of Matrices},
 vol 2, Chelsea Pub. Company, New York, 1964.

 \bibitem{Har} J. Harris, \emph{Algebraic Geometry: A First Course}, Springer, 1992.

 \bibitem{HT1}
 J. Harris and L.W. Tu, On symmetric and skew-symmetric
 determinantal varieties, {\it Topology} {\bf 23} (1984), 71--84.

 \bibitem{Ha70} R Harshman, Foundations of the PARAFAC procedure: Models and conditions for an
 explanatory multi-modal factor analysis, UCLA working papers in phonetics 16 (1970)
 1-84

 \bibitem{Ja79} J. Ja'Ja', Optimal evaluation of a pair of bilinear forms, SIAM J. on Computing 8, (1979)
 281--293.

 \bibitem{KB09} T.G. Kolda and B.W. Bader, Tensor decompositions and applications, \emph{SIAM Review} 51
 (2009) 455--500.

 \bibitem{Kru} J.B. Kruskal, Three-Way Arrays: Rank and Uniqueness
 of Trilinear Decompositions, with Applications to Arithemtic
 Complexity and Statistics, \emph{Linear Algebra Appl.} 18 (1977),
 95-138.

 \bibitem{Kr89} J.B. Kruskal, Rank, decomposition, and uniqueness for 3-way and N-way arrays, in \emph{Multi-
 way data analysis} 7-18, North-Holland, Amsterdam, 1989.

 \bibitem{LMV04} L. de Lathauwer, B. de Moor, and J. Vandewalle,
 Computation of the canonical decomposition by means of a
 simulatenous generalized Schur decomposition, \emph{SIAM J. Matrix
 Anal. Appl.} 26 (2004), 295-327.

 \bibitem{Lic85} T. Lickteig, Typical tensorial rank, \emph{Linear Algebra Appl.}
 69  (1985), 95--120.

 \bibitem{Mum76} D. Mumford, \emph{Algebraic Geometry I: Complex Projective Varieties},
 Springer-Verlag, Berlin, 1976.

 \bibitem{Roc93} R. Rocci, \emph{Unpublished notes}, 1993.

 \bibitem{Sha77} I.R. Shafarevich, \emph{Basic Algebraic Geometry}, Springer-Verlag, Berlin, 1977.

 \bibitem{St06} A. Stegeman, Degeneracy in CANDECOMP/PARAFAC explained for $p \times p \times 2$ arrays of
 rank p + 1 or higher, \emph{Psychometrika} 71 (2006) 483--501.

 \bibitem{Str} V. Strassen,  Rank and optimal computation of generic
 tensors,  \emph{Linear Algebra Appl.}  52/53  (1983), 645--685.

 \bibitem{tB91} J.M.F ten Berge, Kruskal's polynomial for $2\times 2\times 2$ arrays and a generalization $2\times n\times n$ arrays, \emph{Psychometrika} 56 (1991) 631--636

 \bibitem{tB00} J.M.F. Ten Berge, The typical rank of tall three-way arrays, \emph{Psychometrika} 65 (2000) 525-–532.

 \bibitem {tB04} J.M.F. Ten Berge, Partial uniqueness in Candecomp/Parafac, \emph{J. Chemometrics} 18 (2004) 12–-16.

 \bibitem{tBK99} J.M.F. ten Berge and H A L Kiers, Simplicity of core arrays in three-way principal component
 analysis and the typical rank of $p \times q \times 2$ arrays,  \emph{ Linear Algebra Appl.} 294 (1999) 169--179.

 \bibitem{tBS06} J.M.F. ten Berge and A Stegeman, Symmetry transformations for squared sliced three-way
 arrays, with applications to their typical rank, \emph{ Linear Algebra Appl.} 418 (2006) 215--224.


 \bibitem{Ter} A. Terracini, Sulla rappresentazione delle forme
 quaternarie mediante somme di potenze di forme lineari,
 \emph{Atti Rc. Accad. delle Scienze di Torino} 51, 1915-16.

 \bibitem{Tuc} L.R. Tucker, Some mathematical notes of three-mode
 factor analysis, \emph{Psychometrika}, 31 (1966), 279-311.

 \end{thebibliography}
\end{document}